\documentclass[11pt]{article}
\usepackage{amsfonts,amsmath,amssymb,amsthm}
\usepackage[numbers,sort&compress]{natbib}
\usepackage{pgf,tikz,subfigure}
\usetikzlibrary{decorations.markings}

\theoremstyle{plain}
\newtheorem{theorem}{Theorem}[section]
\newtheorem{corollary}[theorem]{Corollary}
\newtheorem{lemma}[theorem]{Lemma}
\newtheorem{proposition}[theorem]{Proposition}
\newtheorem{conjecture}[theorem]{Conjecture}
\theoremstyle{definition}
\newtheorem{definition}[theorem]{Definition}
\theoremstyle{remark}
\newtheorem{remark}[theorem]{Remark}

\begin{document}

\title{$E$-restricted Double Traces}

\author{
Dan Archdeacon\thanks{Deceased (February 2015).} \\
Department of Mathematics \& Statistics, University of Vermont \\
Burlington VT 05405-0156, USA
\and
Luis Goddyn \\
Department of Mathematics, Simon Fraser University \\
Burnaby, B.C. V5A 1S6, Canada \\
goddyn@sfu.ca 
\and 
Jernej Rus
\\
Abelium d.o.o. \\
Jar\v ska cesta 10B, 1000 Ljubljana, Slovenia \\
jernej.rus@abelium.com
}
\date{\today}
\maketitle

\begin{abstract}
For a graph $G$ and $E \subseteq E(G)$, an $E$-restricted strong trace is a closed walk which traverses every edge from $E$ once in each direction and every other edge twice in the same direction, with no nontrivial repetition at any vertex. We characterize graphs admitting $E$-restricted strong traces as well as $E$-restricted $d$-stable traces, as well as directed (mixed graph) versions of these results. We explain how these characterizations upgrade the mathematical model for self-assembling polypeptide nanostructure design first presented by Gradi\v sar et al.\ in 2013.
\end{abstract}


\medskip\noindent
{\bf Math. Subj. Class. (2020):}
05C05, 
05C10, 
05C45, 
92B05, 
92C40, 
94C15. 

\newpage

\section{Introduction}
\label{sec:intro}

In 2013, Gradi\v sar et al.~\cite{gr-2013} presented a novel self-assembly strategy for polypeptide nanostructure design and experimentally demonstrated the formation of a tetrahedron that self-assembles from a single polypeptide chain comprising $12$ concatenated coiled-coil-forming segments separated by flexible peptide hinges. During the construction the path of the polypeptide chain is guided by a prescribed order of segments that traverses each edge of the polyhedron (the tetrahedron in that particular example) exactly twice. The first mathematical model of this process was given by Klav\v zar and Rus~\cite{kl-2013}, who introduced {\em stable traces}: closed walks which traverse every edge of the graph exactly twice, never immediately retrace an edge, and never pass through a vertex twice along the same pair of edges. Fijav\v z et al.~\cite{fi-2014} generalized this model to {\em $d$-stable traces}---closed walks traversing every edge exactly twice such that for every vertex $v$ there is no set $N$ of neighbors of $v$ with $1 \leq |N| \leq d$ for which the walk enters $v$ from $N$ exactly when it exits to a vertex in $N$ (stable traces being the case $d = 2$; the precise convention is fixed in Section~\ref{sec:double})---and to {\em strong traces}, in which no nontrivial set $N$ of any size has this property. Strong traces remedy the two deficiencies of the original model, at vertices of degree at most $2$ and at vertices of degree at least $6$, and connect the problem to single-face embeddings of graphs.

The mathematical model of~\cite{kl-2013, fi-2014} is based on the fact that every polyhedron $P$ composed from a single polymer chain can be naturally represented by the graph $G(P)$ of the polyhedron. Since in the self-assembly process every edge of $G(P)$ corresponds to a coiled-coil dimer, exactly two segments are associated with every edge of $G(P)$. The polyhedron $P$ is then realized by interlocking pairs of polypeptide segments if its graph $G(P)$ contains a closed walk which traverses every edge exactly twice (a double trace). For a polyhedral nanostructure to be stable, that is, not to fall apart or self-assemble into a structure of a different shape than desired, additional conditions are required, and this is why strong traces are used. The two coiled-coil-forming segments on an edge can be aligned either in the same or in opposite directions, which is modeled by parallel and antiparallel edges of a double trace, respectively.

Further use of this self-assembly strategy in~\cite{ko-2015} marked the beginning of a significant resurgence of protein origami, which had spent the better part of the preceding decade overshadowed by DNA origami. This renaissance has since been confirmed: in 2017, Ljubeti\v c et al.~\cite{lj-2017} demonstrated the design and experimental validation of more than twenty distinct polyhedral protein cages---including a tetrahedron, a square pyramid, and a triangular prism---self-assembling from a single polypeptide chain both in vitro and in living cells, establishing coiled-coil protein origami as a mature and versatile platform. While parallel coiled-coil dimers are more frequently characterized for molecular design than antiparallel ones, it has been proven analytically that certain structures---the tetrahedron of~\cite{gr-2013}, for example---cannot be constructed using only parallel or only antiparallel coiled-coil dimers (see Theorems~\ref{thm:anti_strong} and~\ref{thm:parallel_strong} below, both from~\cite{fi-2014}), which necessitates a hybrid model. In this paper we therefore characterize graphs admitting $E$-restricted strong and $d$-stable traces, thereby upgrading the current mathematical model so that it gives fuller control over the self-assembly process and makes the prediction of the properties of the resulting structure easier and more accurate.

Unless stated otherwise, all graphs considered in this paper are connected, finite, simple (without loops and multiple edges), and have at least one edge. We denote the degree of a vertex $v$ by $d_G(v)$, or by $d(v)$ for short if the graph $G$ is clear from the context. The {\em minimum degree} of $G$ is denoted by $\delta(G)$ and the {\em maximum degree} by $\Delta(G)$. If $v$ is a vertex, then $N(v)$ denotes the set of vertices adjacent to $v$ and $E(v)$ the set of edges incident with $v$. Vertices belonging to a subset $V \subseteq V(G)$ are sometimes called {\em $V$-vertices}. A graph in which all vertices have even degree is called an {\em even graph}; analogously, a graph in which all vertices have odd degree is an {\em odd graph}. This should not be confused with the terms {\em even component} and {\em odd component}, which we use for a connected component of a graph with an even, respectively odd, number of edges. A connected even graph is also called an {\em Eulerian graph}, since it admits an Euler tour (a closed walk traversing every edge exactly once). A {\em spanning tree} $T$ of $G$ is a connected acyclic subgraph of $G$ which contains every vertex of $G$. By removing the edges of $T$ from $G$ we obtain the {\em co-tree} $G - E(T)$, which is not necessarily connected; we regard $G - E(T)$ as a spanning subgraph of $G$, so that its isolated vertices form (even) components without edges.

Let $E' \subseteq E(G)$ and let $G[E']$ denote the subgraph of $G$ induced by the edge set $E'$; its vertices are the endvertices of the edges in $E'$. The vertices of $G[E']$ are called {\em $E'$-vertices}; an $E'$-vertex $v$ is a {\em boundary $E'$-vertex} if $v$ is incident with at least one edge from $E'$ and at least one edge from $E(G) \setminus E'$, and an {\em interior $E'$-vertex} otherwise. We denote by $G/E'$ the multigraph obtained from $G$ by replacing every connected component of $G[E']$ with a single vertex while retaining all edges from $E(G) \setminus E'$; more formally, $G/E'$ is obtained by contracting the edges of $E'$ one by one. The new vertices of $G/E'$ are also called {\em $E'$-vertices}; thus the $E'$-vertices of $G/E'$ are precisely the images of the $E'$-vertices of $G$, and a vertex of $G/E'$ which is not an $E'$-vertex is a vertex of $G$ incident with no edge of $E'$. Note that while $G$ is simple, $G/E'$ may contain loops (arising from edges of $E(G) \setminus E'$ whose endvertices lie in the same component of $G[E']$) and multiple edges (arising from edges of $E(G) \setminus E'$ joining the same pair of components of $G[E']$, or the same component and the same vertex). Since $G$ is simple, every loop and every class of multiple edges of $G/E'$ is incident with at least one $E'$-vertex.

For any other terms and concepts from graph theory and topological graph theory not defined here we refer to~\cite{we-1996} and~\cite{moh-2001}, respectively.

\section{Double traces and $E$-restrictions}
\label{sec:double}

A \emph{walk} in a graph $G$ is an alternating sequence 
\begin{equation}
W=w_0 e_1 w_1 \ldots w_{\ell-1} e_\ell w_\ell,
\label{eq:walk}
\end{equation}
such that for each $i=1,\ldots,\ell$, $e_i$ is an edge between the vertices $w_{i-1}$ and $w_i$. When the edges are clear from the context we also write $W = w_0 \ldots w_\ell$. We say that $W$ \emph{passes through} or \emph{traverses} the edges and vertices contained in the sequence~\eqref{eq:walk}. The {\em length} of a walk is the number of edges in the sequence, and $w_0$ and $w_\ell$ are the \emph{endvertices} of $W$. A walk is \emph{closed} if its endvertices coincide.

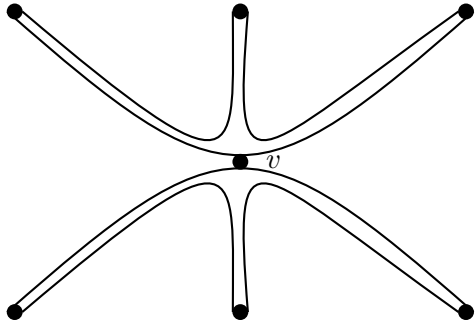
\begin{figure}[ht]
\begin{center}
\begin{tikzpicture}[scale=1.0,style=thick]
\fill (3,2) circle (3pt);
\fill (3.2,2) node[right]{$v$};
\fill (0,0) circle (3pt);
\fill (3,0) circle (3pt);
\fill (6,0) circle (3pt);
\fill (0,4) circle (3pt);
\fill (3,4) circle (3pt);
\fill (6,4) circle (3pt);
\draw (-0.1,4)  .. controls (2.8,1.45) and (3.2,1.45) .. (6.1,4);
\draw (0.1,4)  .. controls (3.1,1.45) and (2.9,2) .. (2.9,4);
\draw (3.1,4)  .. controls (2.9,1.45) and (3.1,2) .. (5.9,4);
\draw (-0.1,0)  .. controls (2.8,2.55) and (3.2,2.55) .. (6.1,0);
\draw (0.1,0)  .. controls (3.1,2.55) and (2.9,2) .. (2.9,0);
\draw (3.1,0)  .. controls (2.9,2.55) and (3.1,2)  .. (5.9,0);
\end{tikzpicture}
\end{center}
\caption{A $3$-repetition at a vertex $v$ of degree $6$: the curves indicate the transitions of the walk at $v$ (after~\cite{fi-2014}).}
\label{fig:repetition}
\end{figure}

A closed walk which traverses every edge of a graph exactly twice is called a {\em double trace}. Applying Euler's theorem to the graph with every edge doubled, many authors observed that every connected graph admits a double trace. Let $W$ be a double trace of length $\ell$ in a graph $G$, $v$ a vertex of $G$ and $N \subseteq N(v)$ a subset of its neighbors. We say that $W$ has an {\em $N$-repetition} at $v$ if at every visit of $W$ to $v$ the walk arrives from a vertex of $N$ exactly when it leaves towards a vertex of $N$; more formally, $W$ has an $N$-repetition at $v$ if
\begin{equation}
\label{eq:repetition}
\text{for every $i \in \{0,\ldots,\ell-1\}$ with $w_i = v$:}\quad w_{i+1} \in N \iff w_{i-1} \in N,
\end{equation}
where indices are taken modulo $\ell$, since we treat a double trace as a cyclic sequence (so $w_1$ is the vertex immediately following $w_\ell = w_0$). An example of a repetition is shown in Figure~\ref{fig:repetition}. An $N$-repetition (at $v$) is a \emph{$d$-repetition} if $|N|=d$, and a $d$-repetition will also be called a repetition \emph{of order $d$}. An $N$-repetition at $v$ is \emph{trivial} if $N=\emptyset$ or $N=N(v)$. Clearly, if $W$ has an $N$-repetition at $v$, then it also has an $(N(v)\setminus N)$-repetition at $v$; we refer to this observation as the \emph{symmetry of repetitions}. In~\cite{fi-2014, rus-2017} a {\em $d$-stable trace} is defined as a double trace without nontrivial repetitions of order at most $d$, and a {\em strong trace} as a double trace without nontrivial repetitions. We keep the definition of a strong trace, but sharpen the definition of a $d$-stable trace in one point: a {\em $d$-stable trace} is a double trace with no $N$-repetition with $1 \leq |N| \leq d$ at all, trivial or not. The two definitions differ only at vertices $v$ with $d(v) \leq d$, where every double trace has the trivial repetition $N(v)$ of order at most $d$: under our definition such vertices exclude $d$-stability, under the literal definition of~\cite{fi-2014, rus-2017} they do not. Our reading is the one implicit in~\cite{fi-2014, rus-2017}, since Proposition~\ref{prp:d-stable} below is stated there as an equivalence with $\delta(G) > d$ (and the remark preceding it there observes that a vertex of degree at most $d$ forces a repetition of order at most $d$ in every double trace); the $\delta(G) > d$ steps in later proofs rely on it. Equivalently, a $d$-stable trace is a double trace without nontrivial repetitions of order at most $d$ in a graph with $\delta(G) > d$. Note that in some other papers the term strong trace is used for (antiparallel) $1$-stable traces; this is not the case in the present paper. It was also observed in~\cite{fi-2014} (for $d = 2$) and~\cite{rus-2017} (for general $d$) that if $\delta(G) > d$, then every strong trace in $G$ is also a $d$-stable trace, and if in addition $\Delta(G) < 2d + 2$, then every $d$-stable trace in $G$ is also a strong trace.

\begin{remark}
\label{rem:Hv}
For a double trace $W$ and a vertex $v$, the {\em local transition multigraph} $H_v$ has vertex set $N(v)$ and, for each visit $\ldots u\, e\, v\, f\, w \ldots$ of $W$ to $v$, one edge $uw$ (loops allowed). Up to replacing neighbors by the corresponding edges at $v$, this is the \emph{vertex figure} $F_{v,W}$ of~\cite{fi-2014}. Two elementary facts are used repeatedly below.
\begin{itemize}
\item[\rm (F1)] A set $N \subseteq N(v)$ is an $N$-repetition of $W$ at $v$ if and only if $N$ is a union of connected components of $H_v$. Hence $W$ has a nontrivial repetition at $v$ if and only if $H_v$ has at least two components, and in that case the minimal nontrivial repetitions at $v$ are the vertex sets of the components of $H_v$, which partition $N(v)$.
\item[\rm (F2)] $H_v$ is $2$-regular (each neighbor $w$ has $H_v$-degree $2$, since the edge $vw$ is traversed exactly twice and each traversal contributes one transition end). Consequently every component of $H_v$ is a cycle (possibly a loop or a pair of multiple edges), and deleting a single edge from a component never disconnects it.
\end{itemize}
If $H_v$ is a single cycle, either of its two orientations is a cyclic permutation of $E(v)$, and if this holds at every vertex, the family of these rotations, together with a suitable edge signature, is a combinatorial embedding of $G$ whose unique facial walk is $W$. This is the connection between strong traces and single-face embeddings on which the proofs in~\cite{fi-2014} rest; in the present paper it enters only through Theorem~\ref{thm:anti_strong}.
\end{remark}

Since every edge $e = uv$ is traversed exactly twice in a double trace $W$, we consider two cases. If $e$ is traversed twice in the same direction (either both times from $u$ to $v$ or both times from $v$ to $u$) then we call $e$ a {\em parallel edge} (with respect to $W$), otherwise $e$ is  an {\em antiparallel edge}. (These notions refer to the traversal by $W$ and should not be confused with multiple edges of a multigraph, which we never call parallel.) A double trace $W$ is a {\em parallel double trace} if every edge of $G$ is parallel and an {\em antiparallel double trace} if every edge of $G$ is antiparallel.

Graphs admitting strong and $d$-stable traces of various kinds were characterized in~\cite{fi-2014, rus-2017} (earlier, $1$-stable and $2$-stable traces were investigated under different names in~\cite{eg-1984, kl-2013, ore-1951, sa-1977, th-1990, tro-1966}), where the following results were proven using a connection between strong traces and single-face embeddings of graphs.

\begin{theorem}
\label{thm:strong}
{\rm \cite[Theorem~2.4]{fi-2014}}
Every connected graph $G$ admits a strong trace.
\end{theorem}

\begin{proposition}
\label{prp:d-stable}
{\rm \cite[Proposition~3.4]{fi-2014}}
Let $G$ be a connected graph. Then $G$ admits a $d$-stable trace if and only if $\delta(G)>d$.
\end{proposition}

\begin{theorem}
\label{thm:anti_strong}
{\rm \cite[Theorem~4.1]{fi-2014}}
A connected graph $G$ admits an antiparallel strong trace if and only if there exists a spanning tree $T$ of $G$ with the property that every connected component of the co-tree $G - E(T)$ has an even number of edges.
\end{theorem}

\begin{theorem}
\label{thm:anti_d-stable}
{\rm \cite[Theorem~1.1]{rus-2017}}
Let $d \geq 1$ be an integer. A connected graph $G$ admits an antiparallel $d$-stable trace if and only if $\delta(G) > d$ and there exists a spanning tree $T$ of $G$ with the property that every connected component of the co-tree $G-E(T)$ is even or contains a vertex $v$ with $d_G(v) \geq 2d + 2$.
\end{theorem}

\begin{theorem}
\label{thm:parallel_strong}
{\rm \cite[Theorem~5.3]{fi-2014}}
A connected graph $G$ admits a parallel strong trace if and only if $G$ is Eulerian.
\end{theorem}

\begin{theorem}
\label{thm:parallel_d-stable}
{\rm \cite[Theorem~5.4]{fi-2014}}
Let $d \geq 1$ be an integer. A connected graph $G$ admits a parallel $d$-stable trace if and only if $G$ is Eulerian and $\delta(G) > d$.
\end{theorem}

We shall also need double traces of multigraphs (loops and multiple edges allowed). Let $G_M$ be a multigraph, $W$ a double trace of $G_M$ and $v$ a vertex of $G_M$. Since a loop at $v$ is entered and left through two distinct ends, repetitions in multigraphs are defined in terms of edge-ends rather than neighbors. Let $\hat{E}(v)$ denote the set of {\em edge-ends} at $v$: every non-loop edge at $v$ contributes one edge-end and every loop at $v$ contributes two, so that $|\hat{E}(v)| = d_{G_M}(v)$. For $N \subseteq \hat{E}(v)$ we say that $W$ has an {\em $N$-repetition} at $v$ if at every visit of $W$ to $v$ the entering edge-end lies in $N$ if and only if the leaving edge-end lies in $N$; this is~\eqref{eq:repetition} with edge-ends in place of neighbors. The notions of trivial and nontrivial repetitions, of the order of a repetition, and of strong and $d$-stable traces then carry over verbatim, and so does Remark~\ref{rem:Hv} with $\hat{E}(v)$ as the vertex set of $H_v$ (each edge-end is used by exactly two transitions, so $H_v$ is again $2$-regular). Likewise, a loop is traversed \emph{from} one of its edge-ends \emph{to} the other, so the notions of parallel and antiparallel edges apply to loops as well: a loop is antiparallel if its two traversals start at different edge-ends. For simple graphs this agrees with the definitions above. Observe also that a vertex $v$ of degree $2$ with neighbors $u$ and $w$ carries no nontrivial repetition in a double trace $W$ if and only if each of the two visits of $W$ to $v$ passes from $u$ to $w$ or from $w$ to $u$ rather than turning back; in this case we say that $W$ \emph{passes straight through} $v$.

The following lemma allows us to move freely between a multigraph and the simple graph obtained by subdividing its loops and multiple edges.

\begin{lemma}
\label{lem:multigraph_strong}
Let $G_M$ be a multigraph and let $G$ be the simple graph constructed from $G_M$ by replacing every loop and every member of every class of multiple edges with an internally vertex-disjoint path of length $3$ and $2$, respectively, and for $v \in V(G_M)$ let the edge-ends of $G_M$ at $v$ be identified with the neighbors of $v$ in $G$ by the natural bijection (a common edge $vw$ with $w$, a multiple edge with the internal vertex of its path, the two ends of a loop with the two internal vertices of its path). Then:
\begin{itemize}
\item[\rm (a)] for every double trace $W_M$ of $G_M$ there is a double trace $W$ of $G$, and
\item[\rm (b)] for every double trace $W$ of $G$ without nontrivial repetitions at the internal vertices of the replacement paths there is a double trace $W_M$ of $G_M$,
\end{itemize}
such that $W$ and $W_M$ traverse common edges in the same directions, $W$ traverses each replacement path in the same directions as $W_M$ traverses the corresponding loop or multiple edge, and for every $v \in V(G_M)$ the local transition multigraphs of $W$ and $W_M$ at $v$ coincide under the natural bijection. In particular, $W$ has a nontrivial repetition of a given order at $v \in V(G_M)$ if and only if $W_M$ does, and $W$ has no nontrivial repetition at the internal vertices of the replacement paths.
\end{lemma}

\begin{proof}
Write the path replacing a loop $q$ at $x$ as $x\, a_q\, b_q\, x$ and the path replacing a multiple edge $f = xy$ as $x\, c_f\, y$.

(a) Given $W_M$, replace each traversal of a loop or multiple edge by the traversal of its replacement path in the corresponding direction (a loop traversed from its first end to its second end becomes $x\, a_q\, b_q\, x$, and in the other direction $x\, b_q\, a_q\, x$). This yields a closed walk $W$ in $G$ which traverses every edge of $G$ exactly twice. Every visit of $W$ to a vertex $v \in V(G_M)$ arises from a visit of $W_M$ to $v$, and its two flanking neighbors are the images under the natural bijection of the two edge-ends used by that visit of $W_M$; hence the local transition multigraphs coincide. At an internal vertex of a replacement path both visits of $W$ pass straight through, so $W$ has no nontrivial repetition there.

(b) Given $W$ without nontrivial repetitions at the internal vertices, $W$ passes straight through every internal vertex, so it traverses every replacement path as a whole, in one direction or the other. Compressing each such traversal into the traversal of the corresponding loop or multiple edge yields a closed walk $W_M$ of $G_M$ which traverses every edge of $G_M$ exactly twice, and the same argument as in (a) shows that the local transition multigraphs at $v \in V(G_M)$ coincide. The remaining assertions follow from (F1).
\end{proof}

We now make the key definition of this paper.

\begin{definition}
\label{def:BR_double}
An {\em $E$-restricted double trace} is a double trace $W$ in a connected graph $G$ such that every edge in $E \subseteq E(G)$ is antiparallel and every edge from $E(G) \setminus E$ is parallel.
\end{definition}

Analogously, for a graph $G$ and $E \subseteq E(G)$, we define an {\em $E$-restricted strong trace} and an {\em $E$-restricted $d$-stable trace} as a strong trace and a $d$-stable trace, respectively, in which every edge from $E$ is antiparallel and every edge from $E(G) \setminus E$ is parallel.

Interestingly, $E$-restricted double traces first appeared more than fifty years ago, when Wagner~\cite{wa-1970} posed (in our language) the problem of characterizing graphs which admit $E$-restricted double traces. The problem was solved by Vestergaard~\cite{ve-1975}; an independent proof appears in Fleischner's monograph~\cite{fl-1991}. The result reads as follows.

\begin{theorem}
\label{thm:BR_double}
{\rm \cite[Theorem~2]{ve-1975}, \cite[Theorem~VIII.13]{fl-1991}}
Let $G$ be a connected graph and $E \subseteq E(G)$. Then $G$ admits an $E$-restricted double trace if and only if $G - E$ is an even graph.
\end{theorem}

In 2017, Wang et al.~\cite{wan-2017, wan2-2017} studied $E$-restricted strong traces under additional assumptions. Under the name $E$-antiparallel strong traces (in~\cite{wan-2017}) and $F$-strong traces (in~\cite{wan2-2017}), they showed that if $E \subsetneq E(G)$ and $E$ is an independent set of edges, or induces a path or a cycle~\cite{wan-2017}, or induces a forest~\cite{wan2-2017}, then $G$ admits an $E$-restricted strong trace if and only if the subgraph $G[E']$, where $E' = E(G) \setminus E$, is an even graph. The hypothesis $E \subsetneq E(G)$ is needed: the boundary case $E = E(G)$ is governed by Theorem~\ref{thm:anti_strong} and is not captured by the even-subgraph condition alone---for a cycle, $G[E']$ is empty and hence vacuously even, yet no antiparallel strong trace exists. Our Theorem~\ref{thm:BR_strong} handles these cases as well; see Corollary~\ref{cor:acyclic_case}. The results of~\cite{wan-2017, wan2-2017} are motivated by the fact that by now more parallel than antiparallel coiled-coil dimers have been characterized for molecular design, so that the case in which $E$ contains relatively few of the edges of $G$ is the most applicable one.

More about double traces in general can be found in~\cite{fl-1990, fl-1991}.

\section{Graphs that admit $E$-restricted strong traces and Xuong tree approximations}
\label{sec:BRdouble_traces}

The main result of this paper reads as follows.

\begin{theorem}
\label{thm:BR_strong}
Let $G$ be a connected graph, $E \subseteq E(G)$, and $E' = E(G) \setminus E$. Then $G$ admits an $E$-restricted strong trace if and only if
\begin{enumerate}
\item[\rm (i)]
the subgraph $G[E']$ is an even graph, and
\item[\rm (ii)]
there exists a spanning tree $T$ of $G/E'$ with the property that every connected component of $G/E' - E(T)$ has an even number of edges or contains an $E'$-vertex.
\end{enumerate}
\end{theorem}

Theorem~\ref{thm:strong} implies that (at least in theory) every connected graph $G$ can be constructed from a single chain of coiled-coil-forming segments. Since the number of coiled-coil-forming segments that can simultaneously interlock into a nanostructure is limited, Theorem~\ref{thm:BR_strong} describes in detail the admissible arrangements of these segments in a chain designed to self-assemble into a desired structure. It also generalizes Theorems~\ref{thm:anti_strong} and~\ref{thm:parallel_strong}, since antiparallel and parallel strong traces are the special cases of $E$-restricted strong traces with $E = E(G)$ and $E = \emptyset$, respectively.

Before presenting the proof of Theorem~\ref{thm:BR_strong} in Section~\ref{sec:proof}, we prove in the rest of this section a few lemmas that are used in that proof.

\begin{lemma}
\label{lemma:joining_cycles}
Let $G$ be a graph and let $W_1$ and $W_2$ be two closed walks in $G$ with a common vertex $v$. Then $W_1$ and $W_2$ can be merged into a single closed walk $W$ which traverses the edges of $W_1$ and $W_2$ with the same multiplicities and in the same directions as $W_1$ and $W_2$, and whose transitions differ from those of $W_1$ and $W_2$ only at $v$.
\end{lemma}

\begin{proof}
Denote the vertices preceding and following $v$ in $W_1$ by $u_1$ and $u_2$, respectively, and the vertices preceding and following $v$ in $W_2$ by $u_3$ and $u_4$, respectively. Let $e_1, e_2, e_3$, and $e_4$ be the edges connecting them to $v$. Let us start walking along $W_1$. When we come to $v$ from $u_1$ along $e_1$, we continue along $e_4$ to $u_4$ and follow $W_2$. When we come back to $v$ from $u_3$ along $e_3$, we continue along $e_2$ to $u_2$ and then follow $W_1$ to its end. This merges the walks $W_1$ and $W_2$ into a single closed walk, see Figure~\ref{fig:construction1}.

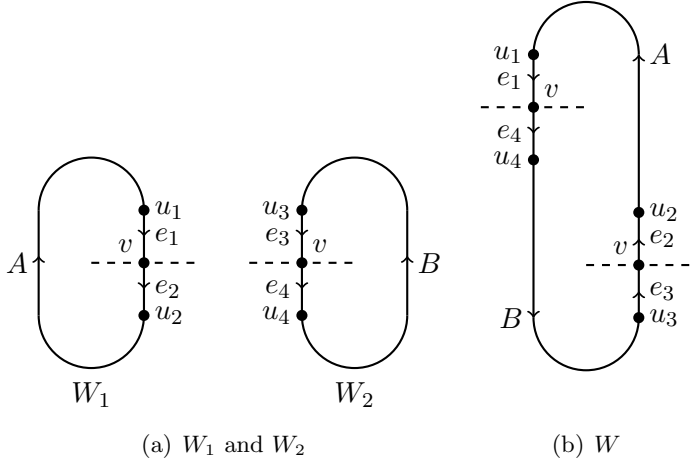
\begin{figure}[ht]
\begin{center}
\subfigure[$W_1$ and $W_2$]
{
\begin{tikzpicture}[scale=0.7,style=thick]
\draw[dashed] (0,2)--(2,2);
\fill (1,3) circle (3pt) node[right]{$u_1$};
\draw (1,2.5) node[right] {$e_1$};
\fill (1,2) circle (3pt) node[above left]{$v$};
\draw (1,1.5) node[right] {$e_2$};
\fill (1,1) circle (3pt) node[right]{$u_2$};
\draw (-1,2) node[left] {$A$};
\draw[decoration={markings, mark=at position 0.5 with {\arrow{>}}},postaction={decorate}] (1,3)--(1,2);
\draw[decoration={markings, mark=at position 0.5 with {\arrow{>}}},postaction={decorate}] (1,2)--(1,1);
\draw[decoration={markings, mark=at position 0.5 with {\arrow{<}}},postaction={decorate}] (-1,3)--(-1,1);
\draw (1,3) arc (0:180:1);
\draw (-1,1) arc (-180:0:1);
\draw (0,-0.5) node {$W_1$};

\draw[dashed] (3,2)--(5,2);
\fill (4,3) circle (3pt) node[left]{$u_3$};
\draw (4,2.5) node[left] {$e_3$};
\fill (4,2) circle (3pt) node[above right]{$v$};
\draw (4,1.5) node[left] {$e_4$};
\fill (4,1) circle (3pt) node[left]{$u_4$};
\draw (6,2) node[right] {$B$};
\draw[decoration={markings, mark=at position 0.5 with {\arrow{>}}},postaction={decorate}] (4,3)--(4,2);
\draw[decoration={markings, mark=at position 0.5 with {\arrow{>}}},postaction={decorate}] (4,2)--(4,1);
\draw[decoration={markings, mark=at position 0.5 with {\arrow{<}}},postaction={decorate}] (6,3)--(6,1);
\draw (6,3) arc (0:180:1);
\draw (4,1) arc (-180:0:1);
\draw (5,-0.5) node {$W_2$};
\end{tikzpicture}
}
\subfigure[$W$]
{
\begin{tikzpicture}[scale=0.7,style=thick]
\draw[dashed] (-2,5)--(0,5);
\draw[dashed] (0,2)--(2,2);
\fill (-1,6) circle (3pt) node[left]{$u_1$};
\draw (-1,5.5) node[left] {$e_1$};
\fill (-1,5) circle (3pt) node[above right]{$v$};
\draw (-1,4.5) node[left] {$e_4$};
\fill (-1,4) circle (3pt) node[left]{$u_4$};

\draw (-1,1) node[left] {$B$};
\fill (1,3) circle (3pt) node[right]{$u_2$};
\draw (1,2.5) node[right] {$e_2$};
\fill (1,2) circle (3pt) node[above left]{$v$};
\draw (1,1.5) node[right] {$e_3$};
\fill (1,1) circle (3pt) node[right]{$u_3$};
\draw (1,6) node[right] {$A$};

\draw[decoration={markings, mark=at position 0.5 with {\arrow{>}}},postaction={decorate}] (-1,6)--(-1,5);
\draw[decoration={markings, mark=at position 0.5 with {\arrow{>}}},postaction={decorate}] (-1,5)--(-1,4);
\draw (-1,4)--(-1,1);
\draw[decoration={markings, mark=at position 0.5 with {\arrow{>}}},postaction={decorate}] (1,1)--(1,2);
\draw[decoration={markings, mark=at position 0.5 with {\arrow{>}}},postaction={decorate}] (1,2)--(1,3);
\draw (1,3)--(1,6);
\draw[decoration={markings, mark=at position 0 with {\arrow{>}}},postaction={decorate}] (-1,1) arc (-180:0:1);
\draw[decoration={markings, mark=at position 0 with {\arrow{>}}},postaction={decorate}] (1,6) arc (0:180:1);
\draw (0,-0.5) node {$\phantom{W}$};
\end{tikzpicture}
}
\end{center}
\caption{Construction of a single closed walk from two closed walks with a common vertex.}
\label{fig:construction1}
\end{figure}
\end{proof}

Note that $v$ need not be the only common vertex of $W_1$ and $W_2$: the merge changes the transitions at $v$ only, and may create a nontrivial repetition there. In Lemma~\ref{lem:assembly} below all merges are performed at $E'$-vertices, where such repetitions are subsequently removed by Corollary~\ref{cor:eliminate}.

The next observation was, to some extent, already presented in~\cite{fi-2014} (see the proof of Theorem~5.3 there).

\begin{lemma}
\label{lemma:vertex_figures_cycles}
Let $W$ be a double trace of a graph $G$ and let $v$ be a vertex at which $W$ has a nontrivial repetition $R_1$. Suppose that some edge $e = uv$ with $u \in R_1$ is traversed by $W$ twice in the same direction. Then there exists a double trace $W'$ of $G$ which traverses every edge in the same directions as $W$ and has strictly fewer nontrivial repetitions at $v$, and hence also strictly fewer in total.
\end{lemma}

\begin{proof}
Set $R_2 := N(v) \setminus R_1$; by the symmetry of repetitions $R_2$ is a nontrivial repetition at $v$, and $R_1 \cap R_2 = \emptyset$. We may assume that $W$ uses the edge $e = uv$, $u \in R_1$, twice in the direction towards $v$; the other case is treated at the end of the proof. Let $e_2=vu_2$ and $e_3=vu_3$ be the edges immediately succeeding the two occurrences of $e$ along $W$ (note that $e_2$ may be equal to $e_3$); since $u \in R_1$ and $R_1$ is a repetition, $u_2, u_3 \in R_1$. Because $R_2 \neq \emptyset$ and every edge to an $R_2$-neighbor is traversed, some visit of $v$ has both flanking neighbors in $R_2$; as the two $e$-visits enter from $u \in R_1$, this visit is distinct from them. Write it $u_4 e_4 v e_5 u_5$ with $e_4 = u_4 v$, $e_5 = vu_5$ and $u_4, u_5 \in R_2$.

Since $W$ is a closed walk, we may cyclically rotate it to begin just after this $R_2$-visit; then the two $e$-traversals and the $R_2$-visit occur in this linear order, and (naming $e_2$ the exit of the first $e$-occurrence met and $e_3$ that of the second---the pair $\{e_2, e_3\}$ being unordered) we may assume that
\[
W = \ldots u\, e\, v\, e_2\, u_2\, A\, u\, e\, v\, e_3\, u_3\, B\, u_4\, e_4\, v\, e_5\, u_5\, C \ldots,
\]
where $A$, $B$, and $C$ are the three ``interior'' subwalks of $W$ between the three shown occurrences of $v$ in $W$. Consider the walk
\[
W' = \ldots u\, e\, v\, e_3\, u_3\, B\, u_4\, e_4\, v\, e_2\, u_2\, A\, u\, e\, v\, e_5\, u_5\, C \ldots
\]
obtained by interchanging the two ``interior'' subwalks $A$ and $B$, see also Figure~\ref{fig:construction2}.

\begin{figure}[ht]
\begin{center}
\subfigure[$W$]
{
\begin{tikzpicture}[scale=1.0,style=thick]
\draw[dashed] (-2,5)--(0,5);
\draw[dashed] (-2,2)--(0,2);
\draw[dashed] (0,3)--(2,3);
\fill (-1,6) circle (3pt) node[left]{$u$};
\draw (-1,5.5) node[left] {$e$};
\fill (-1,5) circle (3pt) node[above right]{$v$};
\draw (-1,4.5) node[left] {$e_2$};
\fill (-1,4) circle (3pt) node[left]{$u_2$};
\draw (-1,3.5) node[left] {$A$};
\fill (-1,3) circle (3pt) node[left]{$u$};
\draw (-1,2.5) node[left] {$e$};
\fill (-1,2) circle (3pt) node[above right]{$v$};
\draw (-1,1.5) node[left] {$e_3$};
\fill (-1,1) circle (3pt) node[left]{$u_3$};

\draw (1,1) node[right] {$B$};
\fill (1,4) circle (3pt) node[right]{$u_5$};
\draw (1,3.5) node[right] {$e_5$};
\fill (1,3) circle (3pt) node[above left]{$v$};
\draw (1,2.5) node[right] {$e_4$};
\fill (1,2) circle (3pt) node[right]{$u_4$};
\draw (1,6) node[right] {$C$};

\draw[decoration={markings, mark=at position 0.5 with {\arrow{>}}},postaction={decorate}] (-1,6)--(-1,5);
\draw[decoration={markings, mark=at position 0.5 with {\arrow{>}}},postaction={decorate}] (-1,5)--(-1,4);
\draw[decoration={markings, mark=at position 0.5 with {\arrow{>}}},postaction={decorate}] (-1,4)--(-1,3);
\draw[decoration={markings, mark=at position 0.5 with {\arrow{>}}},postaction={decorate}] (-1,3)--(-1,2);
\draw[decoration={markings, mark=at position 0.5 with {\arrow{>}}},postaction={decorate}] (-1,2)--(-1,1);
\draw (1,1)--(1,2);
\draw[decoration={markings, mark=at position 0.5 with {\arrow{>}}},postaction={decorate}] (1,2)--(1,3);
\draw[decoration={markings, mark=at position 0.5 with {\arrow{>}}},postaction={decorate}] (1,3)--(1,4);
\draw (1,4)--(1,6);
\draw[decoration={markings, mark=at position 0.75 with {\arrow{>}}},postaction={decorate}] (-1,1) arc (-180:0:1);
\draw[decoration={markings, mark=at position 0 with {\arrow{>}}},postaction={decorate}] (1,6) arc (0:180:1);
\end{tikzpicture}
}
\subfigure[$W'$]
{
\begin{tikzpicture}[scale=1.0,style=thick]
\draw[dashed] (-2,5)--(0,5);
\draw[dashed] (-2,2)--(0,2);
\draw[dashed] (0,3)--(2,3);
\fill (-1,6) circle (3pt) node[left]{$u$};
\draw (-1,5.5) node[left] {$e$};
\fill (-1,5) circle (3pt) node[above right]{$v$};
\draw (-1,4.5) node[left] {$e_3$};
\fill (-1,4) circle (3pt) node[left]{$u_3$};
\draw (-1,3.5) node[left] {$B$};
\fill (-1,3) circle (3pt) node[left]{$u_4$};
\draw (-1,2.5) node[left] {$e_4$};
\fill (-1,2) circle (3pt) node[above right]{$v$};
\draw (-1,1.5) node[left] {$e_2$};
\fill (-1,1) circle (3pt) node[left]{$u_2$};

\draw (1,1) node[right] {$A$};
\fill (1,4) circle (3pt) node[right]{$u_5$};
\draw (1,3.5) node[right] {$e_5$};
\fill (1,3) circle (3pt) node[above left]{$v$};
\draw (1,2.5) node[right] {$e$};
\fill (1,2) circle (3pt) node[right]{$u$};
\draw (1,6) node[right] {$C$};

\draw[decoration={markings, mark=at position 0.5 with {\arrow{>}}},postaction={decorate}] (-1,6)--(-1,5);
\draw[decoration={markings, mark=at position 0.5 with {\arrow{>}}},postaction={decorate}] (-1,5)--(-1,4);
\draw[decoration={markings, mark=at position 0.5 with {\arrow{>}}},postaction={decorate}] (-1,4)--(-1,3);
\draw[decoration={markings, mark=at position 0.5 with {\arrow{>}}},postaction={decorate}] (-1,3)--(-1,2);
\draw[decoration={markings, mark=at position 0.5 with {\arrow{>}}},postaction={decorate}] (-1,2)--(-1,1);
\draw (1,1)--(1,2);
\draw[decoration={markings, mark=at position 0.5 with {\arrow{>}}},postaction={decorate}] (1,2)--(1,3);
\draw[decoration={markings, mark=at position 0.5 with {\arrow{>}}},postaction={decorate}] (1,3)--(1,4);
\draw (1,4)--(1,6);
\draw[decoration={markings, mark=at position 0.75 with {\arrow{>}}},postaction={decorate}] (-1,1) arc (-180:0:1);
\draw[decoration={markings, mark=at position 0 with {\arrow{>}}},postaction={decorate}] (1,6) arc (0:180:1);
\end{tikzpicture}
}
\end{center}
\caption{Decreasing the number of nontrivial repetitions in a double trace (after~\cite{fi-2014}).}
\label{fig:construction2}
\end{figure}
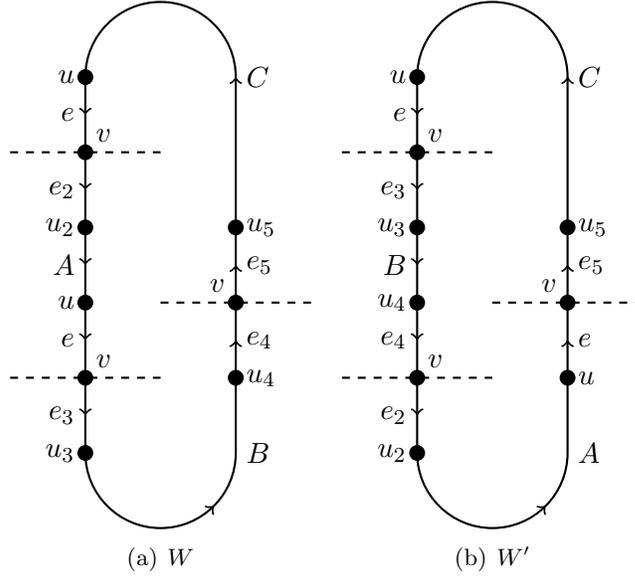

As $W'$ traverses the same collection of edges (in the same directions) as $W$, the walk $W'$ is indeed a double trace which traverses every edge in the same directions as $W$. If $x \ne v$ is a vertex, then the collection of repetitions (trivial or nontrivial) of $W'$ at $x$ equals that of $W$ at $x$, since a pair $e, e'$ of edges meeting at $x$ is consecutive along $W'$ if and only if it is consecutive along $W$.

By construction $W'$ alters only the transitions at $v$, and only those involving $R_1 \cup R_2$. In terms of the local transition multigraph $H_v$ (Remark~\ref{rem:Hv}), the swap deletes the edges $u u_2$ and $u_4 u_5$ and inserts the edges $u u_5$ and $u_2 u_4$. Let $K_1$ be the component of $H_v$ containing $u$ (so $K_1 \subseteq R_1$) and $K_2$ the component containing $u_4$ (so $K_2 \subseteq R_2$); since $R_1 \cap R_2 = \emptyset$, $K_1 \neq K_2$. By (F2) both $K_1$ and $K_2$ are cycles, so $K_1 - u u_2$ and $K_2 - u_4 u_5$ remain connected, and the two new edges $u u_5$, $u_2 u_4$ join them into a single component. Thus $H_v$ for $W'$ has exactly one fewer component than for $W$; by (F1) the number of nontrivial repetitions at $v$ strictly decreases and no new repetition arises (components only merge). By the previous paragraph no other vertex is affected, so the total number of nontrivial repetitions also strictly decreases.

If $e$ is traversed twice in the direction away from $v$, apply the case just proved to the reversed walk $W^{R}$: reversing a closed walk changes neither the local transition multigraph at any vertex nor the set of edges traversed twice in the same direction, and in $W^{R}$ the edge $e$ is traversed twice towards $v$. Reversing the resulting walk once more yields the required $W'$, which traverses every edge in the same directions as $W$.
\end{proof}

The lemma and its proof apply verbatim to double traces of mixed graphs (Section~\ref{sec:directed}) in which every arc is traversed twice in its prescribed direction, arcs being treated like edges; the interchange of subwalks preserves the direction of every traversal.

The hypothesis of Lemma~\ref{lemma:vertex_figures_cycles} that some edge at $v$ is traversed twice in the same direction cannot be dropped. Consider the antiparallel double trace of a cycle which runs once around the cycle and then back; it has nontrivial repetitions only at its starting vertex, namely the two complementary $1$-repetitions there. Since a cycle admits no antiparallel strong trace by Theorem~\ref{thm:anti_strong}, every antiparallel double trace of a cycle has at least two nontrivial repetitions, so no double trace which traverses every edge in the same directions as this one (and which is then again antiparallel) can have fewer nontrivial repetitions than it.

Iterating Lemma~\ref{lemma:vertex_figures_cycles} removes all nontrivial repetitions at any prescribed set of vertices carrying a parallel edge; this is the form in which the lemma is used below.

\begin{corollary}
\label{cor:eliminate}
Let $W$ be a double trace of a graph, or of a mixed graph (see Section~\ref{sec:directed}) in which every arc is traversed twice in its prescribed direction, and let $U$ be a set of vertices each of which is incident with an edge or arc traversed by $W$ twice in the same direction. Then there is a double trace $\widetilde{W}$ which traverses every edge and arc in the same directions as $W$, has no nontrivial repetition at any vertex of $U$, and has the same local transition multigraph as $W$ at every vertex not in $U$ (for mixed graphs, $H_v$ is defined exactly as in Remark~\ref{rem:Hv}, ignoring the directions of arcs).
\end{corollary}

\begin{proof}
Suppose some $v \in U$ has a nontrivial repetition $R$. By the symmetry of repetitions $N(v) \setminus R$ is also a nontrivial repetition, and one of the two contains the neighbor $u$ joined to $v$ by an edge or arc traversed twice in the same direction; take that set as $R_1$ and apply Lemma~\ref{lemma:vertex_figures_cycles}. This preserves all directions and all local transition multigraphs at vertices other than $v$, and strictly decreases the number of nontrivial repetitions at $v$. Induction on the total number of nontrivial repetitions at vertices of $U$ completes the proof.
\end{proof}

The next lemmas extend results about spanning trees presented by Xuong in~\cite{xu-1979, xu2-1979}, where a connection between single-face embeddings of graphs in orientable surfaces and spanning trees was established. For a spanning tree $T$ of a graph $G$, the {\em deficiency} $\xi(G,T)$ of $T$ is the number of odd connected components of the co-tree $G - E(T)$, and $T$ is a {\em Xuong tree} if its deficiency is minimal among all spanning trees of $G$ (this minimum is the {\em deficiency} $\xi(G)$ of the graph $G$). For the purpose of characterizing graphs admitting antiparallel $d$-stable traces, as in Theorem~\ref{thm:anti_d-stable}, one is rather interested in approximations of Xuong trees in which every odd connected component of the co-tree is additionally required to contain a vertex with certain properties. The spanning trees in condition~(ii) of Theorem~\ref{thm:BR_strong} are approximations of Xuong trees in this sense, adapted to $E$-restricted traces.

The following lemma is a multigraph version of~\cite[Lemma~3.2]{rus-2017}.

\begin{lemma}
\label{lem:join_multiple_vertex_trees}
Let $G$ be a graph, $V \subseteq V(G)$, and $T$ a spanning tree of $G$ with the property that every connected component of the co-tree $G - E(T)$ is even or contains a vertex from $V$. Let $V_1, \ldots, V_k$ be pairwise disjoint subsets of $V(G)$ and let $G'$ be the multigraph obtained from $G$ by identifying, for each $i$, the vertices of $V_i$ into a single new vertex $v_i$ (edges are retained with multiplicity; edges inside $V_i$ become loops). Then there exists a spanning tree $T'$ of $G'$ with the property that each connected component of $G' - E(T')$ is even, contains a vertex from $V \setminus (V_1 \cup \cdots \cup V_k)$, or contains one of the new vertices $v_1, \ldots, v_k$.
\end{lemma}

\begin{proof}
Let $\bar{T}$ be the image of $T$ in $G'$ under the identification; it is a connected spanning submultigraph of $G'$ with $|V(G)| - 1$ edges. Every cycle of $\bar{T}$ (loops and pairs of multiple edges included) passes through a new vertex, since otherwise it would lift to a cycle of $T$; remove edges incident with new vertices, one from each cycle, until a spanning tree $T'$ of $G'$ remains. Each removed edge is incident with a new vertex, so it lies in a co-tree component of $G' - E(T')$ containing a new vertex. The identification deletes no edges, so a connected component of $G' - E(T')$ containing no new vertex consists of images of co-tree edges of $G$ between unidentified vertices, and every co-tree edge of $G$ at such a vertex is a co-tree edge of $G'$ with unidentified other end (an identified other end would be a new vertex of the component); hence the component is an isomorphic copy of a connected component of $G - E(T)$, and is even or contains a vertex from $V \setminus (V_1 \cup \cdots \cup V_k)$.
\end{proof}

In this paper Lemma~\ref{lem:join_multiple_vertex_trees} is applied only with $V = \emptyset$ (in the proof of Lemma~\ref{lem:atom_split}), where the sets $V_i$ are the sets of vertices into which a vertex of $G$ was split; the general form costs nothing extra and is recorded for later use.

We next relate spanning trees of a multigraph to spanning trees of the simple graph obtained by subdividing its loops and multiple edges. Recall that a loop is never contained in a spanning tree and that a spanning tree contains at most one edge of each class of multiple edges. The set $X$ in the following lemma is a set of ``exempt'' vertices; in our applications it consists of the $E'$-vertices together, possibly, with all vertices of sufficiently large degree. Note that the degree of a vertex of $G_M$ is the same in $G_M$ and in $G$ below.

\begin{lemma}
\label{lem:multigraph_join_multiple_vertex_trees}
Let $G_M$ be a connected multigraph, $V \subseteq X \subseteq V(G_M)$, and let $G$ be the simple graph obtained from $G_M$ by replacing every loop and every multiple edge with an internally vertex-disjoint path of length $3$ or $2$, respectively. Assume that every loop and every multiple edge of $G_M$ has an endvertex in $V$. Then $G$ has a spanning tree $T$ with the property that every connected component of the co-tree $G - E(T)$ is even or contains a vertex from $X$ if and only if $G_M$ has a spanning tree $T_M$ with the property that every connected component of the co-tree $G_M - E(T_M)$ is even or contains a vertex from $X$.
\end{lemma}

\begin{proof}
Recall that $G$ is obtained from $G_M$ by replacing each loop $q$ at $x$ with a path $x\, a_q\, b_q\, x$ ($3$ edges) and each multiple edge $f = xy$ with a path $x\, c_f\, y$ ($2$ edges), all internally disjoint. We call these paths \emph{gadgets}; the vertices $a_q, b_q, c_f$ are not in $X$.

$(\Rightarrow)$ Let $T$ be as in the statement and $K = G - E(T)$ its co-tree.

\emph{Step 1 (the tree $T_M$).} Define $T_M \subseteq E(G_M)$ as follows: a common edge lies in $T_M$ if and only if it lies in $T$; a multiple edge $f$ lies in $T_M$ if and only if \emph{both} edges of its gadget lie in $T$; loops never lie in $T_M$. Then $T_M$ is a spanning tree of $G_M$: it is acyclic, since a cycle of $T_M$ would lift (replacing each $f \in T_M$ by its gadget) to a cycle of $T$; and it is connected and spanning, since a $T$-path between two vertices of $G_M$ never enters a loop gadget (which meets the rest of $G$ only at $x$) and, whenever it passes through a vertex $c_f$, uses both edges $x c_f, c_f y \in T$, hence projects to $f \in T_M$.

\emph{Step 2 (gadget bookkeeping).} At a loop gadget $x\, a_q\, b_q\, x$, the tree $T$ contains exactly two of the three edges (it must reach $a_q$ and $b_q$, and cannot contain all three). The omitted edge is not $a_q b_q$: otherwise $\{a_q b_q\}$ would be a one-edge (odd) co-tree component containing no vertex from $X$, contradicting the hypothesis. So the omitted edge is $x a_q$ or $b_q x$, a single co-tree edge pendant at $x$, matching the single co-tree edge $\ell$ (a loop at $x$) of $K_M := G_M - E(T_M)$. At a multiple-edge gadget $f = xy$: if both gadget edges lie in $T$, then $f \in T_M$; otherwise exactly one of them does ($T$ must reach $c_f$), and $K$ has one pendant edge through $c_f$ while $K_M$ has $f$ joining $x$ and $y$, so components may merge.

\emph{Step 3 (components of $K_M$).} Let $C_M$ be a connected component of $K_M$. If some vertex $w$ of $C_M$ is an endvertex of a loop or of a class of multiple edges of $G_M$, then a member of that loop or class lies in $K_M$, is incident with $w$, and hence lies in $C_M$; by hypothesis $w$ (or, for a class of multiple edges $xy$, one of $x, y$, both of which then lie in $C_M$) belongs to $V \subseteq X$, so $C_M$ contains a vertex from $X$ and its parity is irrelevant---this is exactly how the component merges of Step 2 are absorbed. Otherwise no vertex of $C_M$ is such an endvertex, so $C_M$ consists only of common edges, present in $K$ with the same endvertices; the component $C$ of $K$ containing them coincides with $C_M$ as an edge set (any $K$-edge at a vertex of $C_M$ would itself be a common edge, hence in $K_M$ and in $C_M$), so $C_M$ is an isomorphic copy of $C$ and is even or contains a vertex from $X$. In all cases $T_M$ has the required property.

$(\Leftarrow)$ Let $T_M$ be as in the statement and $K_M = G_M - E(T_M)$. Let $S$ be the set of all loops and all multiple edges of $G_M$ (every member of every class of multiple edges); these are exactly the edges replaced by gadgets, and each of them has an endvertex in $V$. We insert the gadgets one at a time, in any order; after each insertion we have a multigraph with a spanning tree having the required property, and the single-step argument below uses only that the edge being replaced has an endvertex in $V$. Hence, by induction on $|S|$, it suffices to treat a single edge of $S$.

Let $q$ be a loop at $x \in V$, replaced by the path $x\, a_q\, b_q\, x$. Put $x a_q$ and $a_q b_q$ into the tree and leave $b_q x$ in the co-tree. The result is a spanning tree, and the co-tree changes only in that the loop $q$ (a single co-tree edge at $x$) is replaced by the single pendant co-tree edge $b_q x$ at $x$; every co-tree component keeps its parity and its vertices.

Let $f = xy$ be an edge of $S$ which is not a loop, replaced by the path $x\, c_f\, y$, and say $x \in V$. If $f \in T_M$, put both $x c_f$ and $c_f y$ into the tree; the co-tree is unchanged. If $f \notin T_M$, exactly one of $x c_f$ and $c_f y$ must be put into the tree (to reach $c_f$), and the other one becomes a pendant co-tree edge replacing $f$; all components of $K_M$ other than the component $C_M$ containing $f$ are untouched. If $C_M - f$ is connected, either choice works: the new component has as many edges as $C_M$ and contains all vertices of $C_M$. Otherwise $C_M - f$ consists of a component $C_x$ containing $x$ and a component $C_y$ containing $y$ (either of which may be a single vertex). If $C_y$ is even or contains a vertex from $X$, put $c_f y$ into the tree; the new co-tree components are $C_x + x c_f$, which contains $x \in V \subseteq X$, and $C_y$. If $C_y$ is odd and contains no vertex from $X$, put $x c_f$ into the tree; the new co-tree components are $C_x$, which contains $x \in X$, and $C_y + c_f y$, which is even. In all cases every odd co-tree component contains a vertex from $X$.
\end{proof}

Following~\cite{rus-2017}, \emph{splitting} a vertex $v$ according to a partition $(N_1, \ldots, N_k)$ of $N(v)$ means replacing $v$ by $k$ new pairwise nonadjacent vertices $v_1, \ldots, v_k$ and joining $v_i$ to the vertices of $N_i$ for each $i$; the resulting graph is again simple, and each neighbor of $v$ is adjacent to exactly one of the new vertices.

The following lemma is the common core of the two implications ``trace $\Rightarrow$ tree'' below (in Lemmas~\ref{lemma:tree_with_v_vertices} and~\ref{lemma:tree_with_v_vertices_d}).

\begin{lemma}
\label{lem:atom_split}
Let $W$ be an antiparallel double trace of a connected graph $G$ and let $U$ be the set of vertices at which $W$ has a nontrivial repetition. Then $G$ has a spanning tree $T$ such that every odd connected component of $G - E(T)$ contains a vertex of $U$.
\end{lemma}

\begin{proof}
For $u \in U$ let $N_1^u, \ldots, N_{k_u}^u$ be the vertex sets of the $k_u \geq 2$ components of $H_u$; by (F1) they partition $N(u)$. Obtain a graph $G_U$ from $G$ by replacing every $u \in U$ with $k_u$ new nonadjacent vertices $u_1, \ldots, u_{k_u}$, joining $u_i$ to the vertices of $N_i^u$; each edge formerly incident with $u$ is reattached to exactly one $u_i$. Reroute $W$ accordingly: at every visit of $u$ the two flanking neighbors lie in the same part $N_i^u$ (they are $H_u$-adjacent), so the whole visit passes through $u_i$, and we obtain a closed walk $W_U$ in $G_U$ traversing every edge exactly twice, once in each direction. Hence $G_U$ is connected (and simple, each neighbor of $u$ being adjacent to exactly one $u_i$), every vertex not arising from $U$ keeps its local transition multigraph up to renaming, and $H_{u_i}$ is the connected $i$-th component of $H_u$; so $W_U$ is an antiparallel strong trace of $G_U$, and by Theorem~\ref{thm:anti_strong} $G_U$ has a spanning tree whose co-tree components are all even. Identifying $u_1, \ldots, u_{k_u}$ back into $u$ for every $u \in U$ recovers $G$, and Lemma~\ref{lem:join_multiple_vertex_trees} (with $V = \emptyset$) yields a spanning tree $T$ of $G$ every odd co-tree component of which contains a vertex of $U$.
\end{proof}

The next lemma explains the connection between the above approximations of Xuong trees and antiparallel double traces.

\begin{lemma}
\label{lemma:tree_with_v_vertices}
Let $G$ be a connected graph and $V \subseteq V(G)$. Then $G$ admits an antiparallel double trace whose nontrivial repetitions occur only at vertices of $V$ if and only if there exists a spanning tree $T$ of $G$ with the property that each connected component of $G - E(T)$ is even or contains a vertex from $V$.
\end{lemma}

\begin{proof}
$(\Rightarrow)$ If $W$ is an antiparallel double trace of $G$ whose nontrivial repetitions occur only at vertices of $V$, then Lemma~\ref{lem:atom_split} gives a spanning tree every odd co-tree component of which contains a vertex with a nontrivial repetition, hence a vertex of $V$.

$(\Leftarrow)$ Let $T$ be a spanning tree of $G$ with the property that each connected component of $G - E(T)$ is even or contains a vertex from $V$, and let $\xi = \xi(G,T)$ be its deficiency (the number of odd components of $G - E(T)$; we work with this fixed $T$ and do not require it to be a Xuong tree). We construct, by $\xi$ successive vertex splittings, a graph with a spanning tree of deficiency $0$, and then lift an antiparallel strong trace of that graph back to $G$. If $\xi = 0$ there is nothing to split and Theorem~\ref{thm:anti_strong} applies directly to $G$.

Let $\xi \geq 1$ and let $C$ be an odd connected component of $G - E(T)$; by hypothesis $C$ contains a vertex $v \in V$. Set $E_C := E(v) \cap E(C)$, the co-tree edges at $v$ (all of which lie in $C$), and $E_T := E(v) \setminus E_C$; both are nonempty ($E_T \neq \emptyset$ since $T$ is spanning, $E_C \neq \emptyset$ since $C$ has at least one edge and contains $v$). Construct a new graph $G_1$ from $G$ by replacing $v$ with two new vertices $v'$ and $v''$, attaching the edges of $E_T$ to $v'$ and those of $E_C$ to $v''$; $G_1$ is clearly connected and simple. Let $C_1$ denote $C$ with $v$ renamed $v''$.

For any $e \in E_C$, $T_1 := T + e$ is a spanning tree of $G_1$ (it newly covers $v''$ and creates no cycle); its co-tree is the co-tree of $T$ with $v$ split and $e$ removed, so every component other than $C$ is untouched and $C$ is replaced by $C_1 - e$, which has $|E(C)| - 1$ edges, an even number. We choose $e$ so that all components of $C_1 - e$ are even. If some $e \in E_C$ is not a bridge of $C_1$, choose it: $C_1 - e$ is connected and even. Otherwise every edge at $v''$ is a bridge of $C_1$, so each connected component $D_1, \ldots, D_m$ of $C_1 - v''$ is joined to $v''$ by exactly one edge $e_j$ (two edges to the same $D_j$ would lie on a common cycle). Then $|E(C)| = \sum_{j}(|E(D_j)| + 1)$ is odd, so $|E(D_j)|$ is even for some $j$; choose $e := e_j$. The components of $C_1 - e_j$ are $D_j$, which is even, and the rest, consisting of $v''$ together with the other $D_i$ (connected through $v''$; the isolated vertex $v''$ if $m = 1$), with $|E(C)| - |E(D_j)| - 1$ edges, an even number. In both cases the odd component $C$ is replaced by even components only, so $\xi(G_1, T_1) = \xi(G,T) - 1$ and no new odd component is created; see Figure~\ref{fig:splitting_spanning_tree}. Note that this operation changes the co-tree only inside $C$ and leaves the degrees of all vertices other than $v$ unchanged.

\begin{figure}[ht]
\begin{center}
\subfigure[$G$]
{
\begin{tikzpicture}[scale=1.0,style=thick]
\draw[thin] (-2,1)--(0,0);
\draw[thin] (-2,-1)--(0,0);
\draw[ultra thick] (2,1)--(0,0);
\draw[ultra thick] (2,-1)--(0,0);
\draw[ultra thick] (-2,-1)  .. controls (-1,-2) and (1,-2)  .. (2,-1);
\fill (0,0) circle (3pt);
\fill (0,0.25) node[above]{$v \in V$};
\fill (-2,1) circle (3pt) node[above]{$\notin T$};
\fill (-2,0) node{$\vdots$};
\fill (-2,-1) circle (3pt);
\fill (2,1) circle (3pt) node[above]{$\in T$};
\fill (2,0) node{$\vdots$};
\fill (2,-1) circle (3pt);
\fill (-1,-0.5) node[below]{$e$};
\end{tikzpicture}
}
\subfigure[$G_1$]
{
\begin{tikzpicture}[scale=1.0,style=thick]
\draw[thin] (-2,1)--(-0.5,0);
\draw[ultra thick] (-2,-1)--(-0.5,0);
\draw[ultra thick] (2,1)--(0.5,0);
\draw[ultra thick] (2,-1)--(0.5,0);
\draw[ultra thick] (-2,-1)  .. controls (-1,-2) and (1,-2)  .. (2,-1);
\fill (-0.5,0) circle (3pt);
\fill (-0.75,0) node[left]{$v''$};
\fill (0.5,0) circle (3pt);
\fill (0.75,0) node[right]{$v'$};
\fill (-2,1) circle (3pt);
\fill (-2,0) node{$\vdots$};
\fill (-2,-1) circle (3pt);
\fill (2,1) circle (3pt) node[above]{$\in T_1$};
\fill (2,0) node{$\vdots$};
\fill (2,-1) circle (3pt);
\fill (-1,-0.5) node[below]{$e$};
\end{tikzpicture}
}
\end{center}
\caption{Construction of a spanning tree $T_1$ in $G_1$ from $T$ in $G$. Edges contained in spanning trees are drawn thick.}
\label{fig:splitting_spanning_tree}
\end{figure}
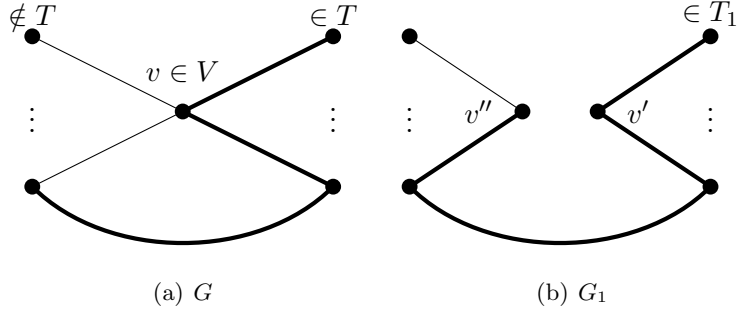

Repeating the same splitting at a $V$-vertex of each of the $\xi$ odd connected components of $G - E(T)$ (each vertex is split at most once, since it lies in at most one component; the splittings do not interfere, since each changes the co-tree only inside its own component), we obtain a graph $G_{\xi}$ with a spanning tree $T_{\xi}$ of deficiency $0$, that is, every connected component of $G_{\xi} - E(T_{\xi})$ has an even number of edges. Therefore, by Theorem~\ref{thm:anti_strong}, $G_{\xi}$ admits an antiparallel strong trace $W_{\xi}$. Let $e = xy$ be an edge of $G_{\xi}$ with $x \notin V(G)$ and let $v \in V$ be the vertex which was replaced by $x$ during the construction; replacing $e$ with $vy$ in $W_{\xi}$, and proceeding in the same way for all edges with endvertices outside $V(G)$, we obtain an antiparallel double trace $W$ of $G$. Finally, its nontrivial repetitions occur only at vertices of $V$: each neighbor $u$ of a split vertex $v$ is adjacent to exactly one of $v', v''$ (the edge $uv$ went to exactly one side of the partition $E_T \sqcup E_C$), so at every unsplit vertex the local transition multigraph is unchanged up to renaming and $W$ has no nontrivial repetition there; at a split vertex $v$, $H_v = H_{v'} \sqcup H_{v''}$ has at most two components, whose only possible nontrivial repetitions are the complementary pair (the $E_T$-side and the $E_C$-side), and $v \in V$. Hence the nontrivial repetitions of $W$ occur only at vertices of $V$.
\end{proof}

Combining Lemmas~\ref{lem:multigraph_strong}, \ref{lem:multigraph_join_multiple_vertex_trees} and~\ref{lemma:tree_with_v_vertices}, we obtain the multigraph version that is actually applied to $G/E'$.

\begin{corollary}
\label{cor:multigraph_tree_trace}
Let $G_M$ be a connected multigraph and $V \subseteq V(G_M)$ such that every loop and every multiple edge of $G_M$ has an endvertex in $V$. Then $G_M$ admits an antiparallel double trace whose nontrivial repetitions occur only at vertices of $V$ if and only if $G_M$ has a spanning tree $T$ such that every connected component of $G_M - E(T)$ is even or contains a vertex from $V$.
\end{corollary}

\begin{proof}
Let $G$ be the simple graph of Lemma~\ref{lem:multigraph_strong}. By that lemma, $G_M$ admits an antiparallel double trace with nontrivial repetitions only at $V$ if and only if $G$ does: part (a) transfers such a trace from $G_M$ to $G$ (the internal path vertices carry no nontrivial repetition), and part (b) transfers it back, since a trace of $G$ with nontrivial repetitions only at $V \subseteq V(G_M)$ has none at the internal path vertices. By Lemma~\ref{lem:multigraph_join_multiple_vertex_trees} with $X = V$, $G_M$ has a spanning tree of the required kind if and only if $G$ does. Lemma~\ref{lemma:tree_with_v_vertices}, applied to $G$, completes the proof.
\end{proof}

In our applications $G_M = G/E'$ (or $G/(E' \cup A)$ in Section~\ref{sec:directed}) and $V$ is the set of $E'$-vertices; the hypothesis of Corollary~\ref{cor:multigraph_tree_trace} is then satisfied because, as observed in Section~\ref{sec:intro}, every loop and every class of multiple edges of $G/E'$ is incident with an $E'$-vertex.

\section{Proof of Theorem~\ref{thm:BR_strong}}
\label{sec:proof}

We first isolate two constructions---projecting an $E$-restricted double trace onto $G/E'$, and lifting an antiparallel double trace of $G/E'$ back to $G$---which are used in the proofs of all four main theorems, and then prove Theorem~\ref{thm:BR_strong}.

\begin{lemma}
\label{lem:projection}
Let $G$ be a connected graph, $E \subseteq E(G)$ nonempty, $E' = E(G) \setminus E$, and let $W$ be a double trace of $G$ in which every edge of $E$ is antiparallel and every edge of $E'$ is parallel. Deleting from $W$ all traversals of edges of $E'$ and mapping every vertex to its image in $G/E'$ yields an antiparallel double trace $W_E$ of $G/E'$. Moreover, for every vertex $v$ of $G/E'$ which is not an $E'$-vertex, the local transition multigraphs of $W_E$ and of $W$ at $v$ coincide (the edge-ends of $G/E'$ at $v$ being the edges of $G$ at $v$).
\end{lemma}

\begin{proof}
The two endvertices of an edge of $E'$ have the same image in $G/E'$, so $W_E$ is a closed walk in $G/E'$. Every edge of $G/E'$ corresponds to a unique edge of $E$, which is traversed once in each direction by $W$; hence $W_E$ is an antiparallel double trace of $G/E'$. A vertex $v$ of $G/E'$ which is not an $E'$-vertex is a vertex of $G$ incident with no edge of $E'$, so no deleted traversal involves $v$: the visits of $W_E$ to $v$ are exactly the visits of $W$ to $v$, and every edge of $G$ at $v$ is retained in $G/E'$ with its end at $v$.
\end{proof}

\begin{lemma}
\label{lem:assembly}
Let $G$ be a connected graph, $E \subseteq E(G)$ and $E' = E(G) \setminus E$ both nonempty, and let $E_1', \ldots, E_k'$ be the connected components of $G[E']$ (we use the same symbol for a component and for its edge set). Let $W^\ast$ be an antiparallel double trace of $G/E'$ and, for each $i$, let $W_i$ be a closed walk in $E_i'$ which visits every vertex of $E_i'$ and traverses every edge of $E_i'$ exactly twice, in the same direction. Then $G$ admits an $E$-restricted double trace $W$ whose local transition multigraph at every vertex which is not an $E'$-vertex coincides with that of $W^\ast$, and at every interior $E'$-vertex of $E_i'$ coincides with that of $W_i$.
\end{lemma}

\begin{proof}
Since $G$ is connected and $E \neq \emptyset$, every component $E_i'$ contains a boundary $E'$-vertex. Every edge of $G/E'$ corresponds to a unique edge of $E$; replacing each edge of $W^\ast$ by the corresponding edge of $G$ we obtain a cyclic sequence of edges of $E$ in which two consecutive edges share an endvertex whenever their common vertex in $G/E'$ is not an $E'$-vertex, while two edges consecutive at an $E'$-vertex of $G/E'$ end at boundary vertices (possibly different ones) of the same component $E_i'$. Cutting this cyclic sequence at every $E'$-vertex (it is cut at least once: $G/E'$ is connected with at least one edge and at least one $E'$-vertex, so some $E'$-vertex is incident with an edge and is visited by $W^\ast$) we obtain a collection $\mathcal{W}_E$ of walks in $G$, called \emph{segments}: each segment starts and ends at a boundary $E'$-vertex, all its internal vertices lie outside $V(G[E'])$, and together the segments traverse every edge of $E$ exactly once in each direction. Moreover, at every vertex $v$ of $G$ which is not an $E'$-vertex, the visits of the segments to $v$ are exactly the visits of $W^\ast$ to $v$.

Let $b$ be a boundary $E'$-vertex. Since $W^\ast$ traverses every edge of $G/E'$ once in each direction, the number of segments ending at $b$ equals the number of segments starting at $b$ (both equal the number of edges of $E$ at $b$). We may therefore pair every segment ending at $b$ with a segment starting at $b$; doing so at every boundary vertex and concatenating paired segments turns $\mathcal{W}_E$ into a collection $Z_1, \ldots, Z_p$ of closed walks in $G$, each passing through at least one boundary $E'$-vertex. Consider the bipartite graph $\mathcal{I}$ with vertex set $\{Z_1, \ldots, Z_p\} \cup \{W_1, \ldots, W_k\}$ in which $Z_j$ and $W_i$ are adjacent if $Z_j$ passes through a vertex of $E_i'$. Such a vertex is a common vertex of $Z_j$ and $W_i$, since $W_i$ visits every vertex of $E_i'$, and it is a boundary $E'$-vertex, since segments meet $V(G[E'])$ only at their endvertices. The graph $\mathcal{I}$ is connected: two segments which are consecutive along $W^\ast$ meet the same $E'$-vertex of $G/E'$, so the closed walks $Z_j$ containing them are both adjacent to the corresponding $W_i$; walking once around $W^\ast$ thus links all the $Z_j$, and every $W_i$ is adjacent to some $Z_j$ because $E_i'$ has a boundary vertex. Choose a spanning tree of $\mathcal{I}$ and merge the walks along its edges, one at a time, by Lemma~\ref{lemma:joining_cycles}, each time at a common boundary $E'$-vertex of the two walks being merged. (Two walks $Z_j$, $Z_{j'}$ are never merged directly, as they may share a vertex which is not an $E'$-vertex.) The result is a single closed walk $W$ in $G$ which traverses every edge of $E'$ twice in the same direction (as the $W_i$ do) and every edge of $E$ once in each direction (as the segments do); that is, $W$ is an $E$-restricted double trace of $G$. Since all merges take place at boundary $E'$-vertices, the local transition multigraph of $W$ at a vertex which is not an $E'$-vertex is that of $W^\ast$, and at an interior $E'$-vertex of $E_i'$ it is that of $W_i$.
\end{proof}

Lemmas~\ref{lem:projection} and~\ref{lem:assembly} apply verbatim to a simple mixed graph $G$ with $E' \cup A(G)$ in place of $E'$, where every arc is traversed twice in its prescribed direction and hence behaves like an edge of $E'$; in Lemma~\ref{lem:assembly} the walks $W_i$ then traverse every arc of $E_i'$ twice in its prescribed direction, and so does $W$.

\begin{proof}[Proof of Theorem~\ref{thm:BR_strong}]
$(\Rightarrow)$ Let $W$ be an $E$-restricted strong trace of $G$. Since $W$ is in particular an $E$-restricted double trace, $G - E$ is an even graph by Theorem~\ref{thm:BR_double}; since $G - E$ is $G[E']$ together with isolated vertices, this is (i). If $E = \emptyset$, then $G/E'$ consists of a single vertex without edges and (ii) holds trivially. Otherwise Lemma~\ref{lem:projection} yields an antiparallel double trace $W_E$ of $G/E'$ whose local transition multigraph at every vertex which is not an $E'$-vertex is that of $W$; as $W$ is a strong trace, the nontrivial repetitions of $W_E$ occur only at $E'$-vertices, and Corollary~\ref{cor:multigraph_tree_trace} gives (ii).

$(\Leftarrow)$ Assume that (i) and (ii) hold. If $E' = \emptyset$, then $G/E' = G$ and Corollary~\ref{cor:multigraph_tree_trace} with $V = \emptyset$ gives an antiparallel strong trace of $G$, which is an $E$-restricted strong trace for $E = E(G)$. If $E = \emptyset$, then $G$ is Eulerian by (i) and admits a parallel strong trace by Theorem~\ref{thm:parallel_strong}. Hence assume that $E \neq \emptyset$ and $E' \neq \emptyset$, and let $E_1', \ldots, E_k'$ be the connected components of $G[E']$. Each $E_i'$ is a connected even graph, so it has an Euler tour; traversing it twice gives a closed walk $W_i$ which visits every vertex of $E_i'$ and traverses every edge of $E_i'$ twice in the same direction. By (ii) and Corollary~\ref{cor:multigraph_tree_trace}, $G/E'$ admits an antiparallel double trace $W^\ast$ whose nontrivial repetitions occur only at $E'$-vertices. Lemma~\ref{lem:assembly} yields an $E$-restricted double trace $W$ of $G$ whose local transition multigraph at every vertex which is not an $E'$-vertex is that of $W^\ast$; hence $W$ has no nontrivial repetition at such vertices. Finally, every $E'$-vertex is incident with an edge of $E'$, which is traversed twice in the same direction, so Corollary~\ref{cor:eliminate} with $U$ the set of $E'$-vertices turns $W$ into an $E$-restricted double trace without nontrivial repetitions, that is, into an $E$-restricted strong trace.
\end{proof}

\begin{corollary}
\label{cor:acyclic_case}
Let $G$ be a connected graph, $E \subseteq E(G)$ and $E' = E(G) \setminus E$, and suppose that every cycle of $G[E]$ contains an $E'$-vertex. Then $G$ admits an $E$-restricted strong trace if and only if $G[E']$ is an even graph.
\end{corollary}

\begin{proof}
By Theorem~\ref{thm:BR_strong} it suffices to show that condition~(ii) holds automatically. The edges of $E$ joining two vertices which are not $E'$-vertices are untouched by the contraction $G/E'$ and form an acyclic subgraph $F_0$ of $G/E'$, since a cycle in $F_0$ would be a cycle of $G[E]$ without an $E'$-vertex. Extend $F_0$ to a spanning tree $T$ of $G/E'$. Every edge of $G/E' - E(T)$ is then incident with an $E'$-vertex, so every connected component of $G/E' - E(T)$ with at least one edge contains an $E'$-vertex, and the components without edges are even.
\end{proof}

The hypothesis of Corollary~\ref{cor:acyclic_case} is satisfied whenever $G[E]$ is a forest, and also whenever $G[E]$ is a cycle and $E \neq E(G)$ (a cycle of $G[E]$ avoiding all $E'$-vertices consists of vertices of degree $2$ in $G$ and is therefore a connected component of $G$, so $G$ itself is that cycle and $E = E(G)$). Corollary~\ref{cor:acyclic_case} thus contains the characterizations of~\cite{wan-2017} ($E$ an independent set of edges, or inducing a path or a cycle) and of~\cite{wan2-2017} ($E$ inducing a forest). The hypothesis cannot be dropped: for the graph obtained from two disjoint triangles $abc$ and $xyz$ by adding the edge $cx$, with $E' = \{xy, yz, zx\}$, the subgraph $G[E']$ is even, but every spanning tree of $G/E'$ leaves a single edge of the triangle $abc$ in the co-tree, an odd component without an $E'$-vertex, so condition~(ii) of Theorem~\ref{thm:BR_strong} fails.

\section{Generalization to $d$-stable traces}
\label{sec:dstable}

\begin{theorem}
\label{thm:BR_d-stable}
Let $G$ be a connected graph, $d$ a positive integer, $E \subseteq E(G)$, and $E' = E(G) \setminus E$. Then $G$ admits an $E$-restricted $d$-stable trace if and only if $\delta(G) > d$ and
\begin{enumerate}
\item[\rm (i)]
the subgraph $G[E']$ is an even graph, and
\item[\rm (ii)]
there exists a spanning tree $T$ of $G/E'$ with the property that every connected component of $G/E' - E(T)$ is even, or contains a vertex $v$ with $d_{G/E'}(v) \geq 2d + 2$, or contains an $E'$-vertex.
\end{enumerate}
\end{theorem}

Note that a vertex $v$ of $G/E'$ which is not an $E'$-vertex satisfies $d_{G/E'}(v) = d_G(v)$; since $E'$-vertices are exempt in (ii) anyway, the degree condition there may equivalently be read in $G$.

Before proving Theorem~\ref{thm:BR_d-stable} we establish the $d$-stable analogue of Lemma~\ref{lemma:tree_with_v_vertices}. Its proof cannot simply invoke Theorem~\ref{thm:anti_d-stable}, because that theorem requires $\delta > d$, which fails for the vertices created by the splittings in the proof of Lemma~\ref{lemma:tree_with_v_vertices}. Instead we use the following balanced splitting lemma from~\cite{rus-2017}, where splitting is as defined before Lemma~\ref{lem:atom_split}.

\begin{lemma}[{\cite[Lemma~3.4]{rus-2017}}, stated there with $d$ for $t$ and $\xi(G,d,T)$ for $\xi(G,T)$]
\label{lem:balanced_split}
Let $G$ be a connected graph, $t \geq 1$ an integer, and $T$ a spanning tree of $G$ such that every connected component of $G - E(T)$ is even or contains a vertex of degree at least $t$. Let $v$ be a vertex of degree at least $t$ contained in an odd component of $G - E(T)$. Then there is a partition of $N(v)$ into two sets of sizes $\lceil d(v)/2 \rceil$ and $\lfloor d(v)/2 \rfloor$ such that the graph $G_1$ obtained from $G$ by splitting $v$ according to this partition is connected and has a spanning tree $T_1$ such that every connected component of $G_1 - E(T_1)$ is even or contains a vertex of degree at least $t$, and $\xi(G_1, T_1) < \xi(G, T)$.
\end{lemma}

In~\cite{rus-2017} this lemma is derived from Lemma~3.3 there, which provides the balanced splitting together with a spanning tree $T_1$ satisfying $\xi(G_1, T_1) < \xi(G, T)$. That every odd component of $G_1 - E(T_1)$ still contains a vertex of degree at least $t$ can be read off from that proof as follows. In each of its cases the co-tree is modified only by moving edges incident with $v$ into the tree, which cuts the odd component containing $v$ into pieces, and by moving at most two tree edges sharing an endvertex into the co-tree, which merges some of these pieces with components of $G - E(T)$; in the last subcase this happens in two rounds, via an auxiliary spanning tree of $G$ with three tree edges at $v$. All pieces cut off from the component of $v$ are even, except possibly the piece containing $v$ itself, which is cut into even pieces in the final round, and a merged component has the parity of the sum of its constituents, since exactly two edges are added to it. Hence every odd component of the final co-tree contains an odd component of $G - E(T)$, and with it a vertex of degree at least $t$ whose degree is not affected by the splitting of $v$.

\begin{lemma}
\label{lemma:tree_with_v_vertices_d}
Let $G$ be a connected graph, $V \subseteq V(G)$, and $d \geq 1$ an integer. Then $G$ admits an antiparallel double trace whose nontrivial repetitions of order at most $d$ occur only at vertices of $V$ if and only if there exists a spanning tree $T$ of $G$ such that every connected component of $G - E(T)$ is even, or contains a vertex of $V$, or contains a vertex $w$ with $d_G(w) \geq 2d + 2$.
\end{lemma}

\begin{proof}
$(\Rightarrow)$ Let $W$ be such an antiparallel double trace and $U$ the set of vertices at which $W$ has a nontrivial repetition. By Lemma~\ref{lem:atom_split}, $G$ has a spanning tree every odd co-tree component of which contains a vertex $u \in U$. If $u \notin V$, then every nontrivial repetition at $u$ has order at least $d+1$; the vertex sets of the components of $H_u$ are nontrivial repetitions by (F1), there are at least two of them, and they partition $N(u)$, so $d_G(u) \geq 2(d+1) = 2d + 2$.

$(\Leftarrow)$ Suppose that such a spanning tree $T$ exists. We eliminate the odd co-tree components in two stages, splitting vertices as we go; throughout, the current graph is simple and connected. In the first stage, for every odd component $C$ of $G - E(T)$ which contains a vertex $v \in V$, split $v$ by the bridge argument in the proof of Lemma~\ref{lemma:tree_with_v_vertices}. As noted there, this replaces $C$ by even components, changes the co-tree only inside $C$, and leaves the degrees of all other vertices unchanged; since a vertex lies in at most one co-tree component, every vertex is split at most once in this stage. After the first stage, every odd component of the current co-tree contains a vertex of degree at least $2d + 2$ (its exempt vertex was not in $V$ and hence was not split). In the second stage, apply Lemma~\ref{lem:balanced_split} with $t = 2d + 2$ repeatedly, each time at a vertex of degree at least $2d + 2$ in an odd component. The conclusion of the lemma is again its hypothesis and the deficiency strictly decreases, so after finitely many steps we arrive at a connected graph $G_\xi$ with a spanning tree whose co-tree components are all even. (In the second stage a vertex may be split several times; each time a vertex of degree at least $2d + 2$ is split into two vertices of degree at least $d + 1$.) By Theorem~\ref{thm:anti_strong}, $G_\xi$ admits an antiparallel strong trace $W_\xi$. Identifying, for every vertex $u$ of $G$, all vertices of $G_\xi$ that arose from $u$ back into $u$, as in the proof of Lemma~\ref{lemma:tree_with_v_vertices}, turns $W_\xi$ into an antiparallel double trace $W$ of $G$. A vertex of $G$ that was never split has the same local transition multigraph in $W$ as in $W_\xi$ and hence carries no nontrivial repetition. For a split vertex $u$, the local transition multigraph $H_u$ of $W$ is the disjoint union of the (connected) local transition multigraphs of the vertices of $G_\xi$ that arose from $u$, so by (F1) every nontrivial repetition of $W$ at $u$ is a union of the neighborhoods of some of these vertices. If $u \in V$, repetitions of any order are permitted at $u$. If $u \notin V$, then $u$ was split only in the second stage, so every vertex of $G_\xi$ arising from $u$ has degree at least $d + 1$, and consequently every nontrivial repetition of $W$ at $u$ has order at least $d + 1$. Therefore $W$ has no nontrivial repetition of order at most $d$ outside $V$.
\end{proof}

Exactly as Corollary~\ref{cor:multigraph_tree_trace} was derived from Lemma~\ref{lemma:tree_with_v_vertices}, Lemma~\ref{lemma:tree_with_v_vertices_d} yields its multigraph version; here one applies Lemma~\ref{lem:multigraph_join_multiple_vertex_trees} with $X$ the set of vertices of $V$ together with all vertices of degree at least $2d + 2$, which is legitimate because the internal vertices of the replacement paths have degree $2 < 2d+2$ and all other degrees agree in $G_M$ and $G$. Lemma~\ref{lem:multigraph_strong}(b) is applicable because a nontrivial repetition at an internal vertex of degree $2$ has order $1 \leq d$, and such repetitions are excluded outside $V$.

\begin{corollary}
\label{cor:multigraph_tree_trace_d}
Let $G_M$ be a connected multigraph, $d \geq 1$ an integer, and $V \subseteq V(G_M)$ such that every loop and every multiple edge of $G_M$ has an endvertex in $V$. Then $G_M$ admits an antiparallel double trace whose nontrivial repetitions of order at most $d$ occur only at vertices of $V$ if and only if $G_M$ has a spanning tree $T$ such that every connected component of $G_M - E(T)$ is even, or contains a vertex from $V$, or contains a vertex $w$ with $d_{G_M}(w) \geq 2d + 2$.
\end{corollary}

We remark that Corollary~\ref{cor:multigraph_tree_trace_d} with $V = \emptyset$, together with the observation that under $\delta(G) > d$ no trivial repetition $N(v)$ with $d(v) \leq d$ can occur, yields an alternative derivation of Theorem~\ref{thm:anti_d-stable} from Lemma~\ref{lem:balanced_split}.

\begin{proof}[Proof of Theorem~\ref{thm:BR_d-stable}]
$(\Rightarrow)$ Let $W$ be an $E$-restricted $d$-stable trace of $G$. Condition (i) follows from Theorem~\ref{thm:BR_double}, and $\delta(G) > d$ from Proposition~\ref{prp:d-stable} (recall that the trivial repetition $N(v)$ is forbidden at vertices with $d(v) \leq d$). If $E = \emptyset$, then (ii) holds trivially. Otherwise Lemma~\ref{lem:projection} yields an antiparallel double trace $W_E$ of $G/E'$ whose local transition multigraph at every vertex which is not an $E'$-vertex is that of $W$; as $W$ is $d$-stable, the nontrivial repetitions of $W_E$ of order at most $d$ occur only at $E'$-vertices, and Corollary~\ref{cor:multigraph_tree_trace_d} gives (ii).

$(\Leftarrow)$ Assume that $\delta(G) > d$ and that (i) and (ii) hold. If $E' = \emptyset$, then $G/E' = G$ and Corollary~\ref{cor:multigraph_tree_trace_d} with $V = \emptyset$ gives an antiparallel double trace of $G$ without nontrivial repetitions of order at most $d$; since $\delta(G) > d$, it has no trivial repetition $N(v)$ with $d(v) \leq d$ either, so it is an antiparallel $d$-stable trace. If $E = \emptyset$, then $G$ is Eulerian and admits a parallel $d$-stable trace by Theorem~\ref{thm:parallel_d-stable}. Hence assume that $E \neq \emptyset$ and $E' \neq \emptyset$, and let $W_1, \ldots, W_k$ be doubled Euler tours of the components of $G[E']$ as in the proof of Theorem~\ref{thm:BR_strong}. By (ii) and Corollary~\ref{cor:multigraph_tree_trace_d}, $G/E'$ admits an antiparallel double trace $W^\ast$ whose nontrivial repetitions of order at most $d$ occur only at $E'$-vertices. Lemma~\ref{lem:assembly} yields an $E$-restricted double trace of $G$ whose local transition multigraph at every vertex which is not an $E'$-vertex is that of $W^\ast$, and Corollary~\ref{cor:eliminate} with $U$ the set of $E'$-vertices turns it into an $E$-restricted double trace $\widetilde{W}$ without nontrivial repetitions at $E'$-vertices and with the local transition multigraphs of $W^\ast$ elsewhere. Thus every nontrivial repetition of $\widetilde{W}$ has order greater than $d$, and since $\delta(G) > d$ there is no trivial repetition $N(v)$ with $d(v) \leq d$; that is, $\widetilde{W}$ is an $E$-restricted $d$-stable trace.
\end{proof}

\section{Directed versions of Theorems~\ref{thm:BR_strong} and~\ref{thm:BR_d-stable}}
\label{sec:directed}

In this section we prove directed versions of Theorems~\ref{thm:BR_strong} and~\ref{thm:BR_d-stable}. We first recall some definitions concerning mixed graphs. A {\em mixed graph} $G$ consists of a set of vertices $V(G)$, a set of (undirected) edges $E(G)$, and a set of directed edges or {\em arcs} $A(G)$. All mixed graphs considered here are finite and simple, in the sense that there are no loops and any two vertices are joined by at most one edge or arc. For every $X \subseteq V(G)$, the {\em cut quantities} $e(X)$, $a^+(X)$, and $a^-(X)$ denote, respectively, the number of edges with exactly one endvertex in $X$, the number of arcs directed from $X$ to $V(G) \setminus X$, and the number of arcs directed from $V(G) \setminus X$ to $X$. The {\em degree} $d_G(v)$ of a vertex counts all incident edges and arcs, and $\delta(G)$ is the minimum degree; for a singleton $X = \{v\}$ we write $e(v)$, $a^+(v)$, $a^-(v)$ for the numbers of edges, out-arcs, and in-arcs at $v$. A mixed graph is {\em weakly connected} if the graph obtained by replacing all arcs with undirected edges is connected, and the {\em connected components} of a mixed graph are those of this underlying graph. For a set $S \subseteq E(G) \cup A(G)$, $G_S$ denotes the sub-mixed-graph with edge and arc set $S$ whose vertices are the endvertices of the elements of $S$, and $G/S$ denotes the multigraph obtained from $G$ by contracting every connected component of $G_S$ to a single vertex while retaining the edges of $E(G) \setminus S$. As in Section~\ref{sec:intro}, the vertices of $G_S$ and the new vertices of $G/S$ are called {\em $S$-vertices}, a vertex of $G_S$ incident with an edge of $E(G) \setminus S$ is a {\em boundary $S$-vertex}, and since $G$ is simple every loop and every class of multiple edges of $G/S$ is incident with an $S$-vertex. A walk in a mixed graph may traverse an edge in either direction but an arc only in its prescribed direction. An {\em Euler tour} of a mixed graph is a closed walk traversing every edge and every arc exactly once, and a {\em cycle decomposition} is a partition of $E(G) \cup A(G)$ into cycles each of which can be traversed as a closed walk in this sense. Double traces and their repetitions, as well as parallel and antiparallel edges, are defined for mixed graphs exactly as for graphs, a double trace now traversing every edge and every arc exactly twice; we consider only double traces which traverse every arc twice in its prescribed direction, so that arcs behave like parallel edges. In particular, for $E \subseteq E(G)$ and $E' = E(G) \setminus E$, an {\em $E$-restricted strong trace} (resp.\ {\em $E$-restricted $d$-stable trace}) of a mixed graph $G$ is a strong trace (resp.\ $d$-stable trace) of $G$ in which every edge of $E$ is antiparallel, every edge of $E'$ is parallel, and every arc is traversed twice in its prescribed direction.

The following theorem was proven by Ford and Fulkerson~\cite{fo-1962} and rediscovered by Batagelj and Pisanski~\cite{ba-1979} and by Fleischner~\cite{fl-1983}.

\begin{theorem}
\label{thm:max_flow}
{\rm \cite[Theorem~7.1]{fo-1962}, \cite{ba-1979}, \cite[Theorem~IV.11]{fl-1983}}
Let $G$ be a weakly connected mixed graph. Then the following three statements are equivalent:
\begin{enumerate}
\item
For every $X \subseteq V(G)$, $f_G(X) = e(X)- |a^+(X) - a^-(X)|$ is a non-negative even integer.
\item
$G$ has an Euler tour.
\item
$G$ has a cycle decomposition.
\end{enumerate}
\end{theorem}

We first characterize mixed graphs admitting $E$-restricted strong traces.

\begin{theorem}
\label{thm:BR_directed_strong}
Let $G$ be a weakly connected mixed graph with arc set $A = A(G)$, $E \subseteq E(G)$, and $E' = E(G) \setminus E$. Then $G$ admits an $E$-restricted strong trace (in which every arc is traversed twice in its prescribed direction) if and only if
\begin{enumerate}
\item[\rm (i)]
for every connected component $B$ of $G_{E' \cup A}$ and every $X \subseteq V(B)$, $f_{B}(X) = e_B(X) - |a_B^+(X) - a_B^-(X)|$ is a non-negative even integer, and
\item[\rm (ii)]
there exists a spanning tree $T$ of $G/(E' \cup A)$ with the property that every connected component of $G/(E' \cup A) - E(T)$ has an even number of edges or contains an $(E' \cup A)$-vertex.
\end{enumerate}
\end{theorem}

\begin{remark}
\label{rem:mixed_vertex_condition}
Condition~(i) of Theorem~\ref{thm:BR_directed_strong} (and of Theorem~\ref{thm:BR_directed_d-stable} below) applied to the singletons $X = \{v\}$ gives the following explicit vertex-wise necessary condition: for every vertex $v$ of $G_{E' \cup A}$,
\[
e(v) \geq |a^+(v) - a^-(v)| \quad \text{and} \quad e(v) \equiv |a^+(v) - a^-(v)| \pmod{2},
\]
where the cut quantities are taken in the component of $G_{E' \cup A}$ containing $v$. This is the mixed-graph analogue of the requirement that $G[E']$ be an even graph: in the undirected case, even degree at each vertex of $G[E']$ ensures that the edges of $E'$ can be oriented so that every vertex is balanced, while in the mixed case the $e(v)$ edges of $E'$ at $v$ must be able to compensate the imbalance $|a^+(v) - a^-(v)|$ of the arcs at $v$---the inequality ensures that there are enough of them, and the parity condition ensures that the remaining ones can be oriented in pairs. We highlight this special case since it is the most immediately checkable consequence of (i).
\end{remark}

We prove Theorem~\ref{thm:BR_directed_strong} together with its $d$-stable version, which we state first.

\begin{theorem}
\label{thm:BR_directed_d-stable}
Let $G$ be a weakly connected mixed graph with arc set $A = A(G)$, $d$ a positive integer, $E \subseteq E(G)$, and $E' = E(G) \setminus E$. Then $G$ admits an $E$-restricted $d$-stable trace (in which every arc is traversed twice in its prescribed direction) if and only if $\delta(G) > d$ and
\begin{enumerate}
\item[\rm (i)]
for every connected component $B$ of $G_{E' \cup A}$ and every $X \subseteq V(B)$, $f_{B}(X) = e_B(X) - |a_B^+(X) - a_B^-(X)|$ is a non-negative even integer, and
\item[\rm (ii)]
there exists a spanning tree $T$ of $G/(E' \cup A)$ with the property that every connected component of $G/(E' \cup A) - E(T)$ is even, or contains a vertex $v$ with $d_{G/(E' \cup A)}(v) \geq 2d + 2$, or contains an $(E' \cup A)$-vertex.
\end{enumerate}
\end{theorem}

\begin{proof}[Proof of Theorems~\ref{thm:BR_directed_strong} and~\ref{thm:BR_directed_d-stable}]
Write $S = E' \cup A$; every loop and every class of multiple edges of $G/S$ is incident with an $S$-vertex, so Corollaries~\ref{cor:multigraph_tree_trace} and~\ref{cor:multigraph_tree_trace_d} apply to $G/S$ with $V$ the set of $S$-vertices.

$(\Rightarrow)$ Let $W$ be an $E$-restricted strong trace (resp.\ $d$-stable trace) of $G$ traversing every arc twice in its prescribed direction. Fix a connected component $B$ of $G_S$. At each vertex $v \in V(B)$, the edges of $E$ incident with $v$ are antiparallel and therefore contribute equally to the entries into and exits from $v$; as $W$ enters and leaves $v$ equally often, the traversals of edges and arcs of $B$ are balanced at $v$ as well. Hence the orientation of $B$ induced by $W$ (each edge of $E'$ oriented in its common direction of traversal, each arc in its prescribed direction) has equal in- and out-degree at every vertex, because each edge and arc of $B$ is traversed exactly twice. A weakly connected mixed graph with a balanced orientation has an Euler tour, so Theorem~\ref{thm:max_flow}, (2)$\Rightarrow$(1), applied to $B$ yields (i). In the $d$-stable case, $\delta(G) > d$ follows from Proposition~\ref{prp:d-stable} applied to the underlying graph of $G$, of which $W$ is a $d$-stable trace. For (ii), we may assume $E \neq \emptyset$ (otherwise $G/S$ is a single vertex). Lemma~\ref{lem:projection}, applied with $S$ in place of $E'$, yields an antiparallel double trace of $G/S$ whose local transition multigraph at every vertex which is not an $S$-vertex is that of $W$; hence its nontrivial repetitions (resp.\ its nontrivial repetitions of order at most $d$) occur only at $S$-vertices, and Corollary~\ref{cor:multigraph_tree_trace} (resp.\ Corollary~\ref{cor:multigraph_tree_trace_d}) gives (ii).

$(\Leftarrow)$ If $S = \emptyset$, then $G$ is a graph and $G/S = G$, so Corollary~\ref{cor:multigraph_tree_trace} (resp.\ Corollary~\ref{cor:multigraph_tree_trace_d}) with $V = \emptyset$ gives an antiparallel strong trace (resp.\ an antiparallel double trace without nontrivial repetitions of order at most $d$, which is $d$-stable since $\delta(G) > d$ excludes trivial repetitions $N(v)$ with $d(v) \leq d$). If $E = \emptyset$, then $G_S = G$ is weakly connected and, by (i) and Theorem~\ref{thm:max_flow}, has an Euler tour; traversing it twice gives a double trace in which every edge and arc is traversed twice in the same direction, and Corollary~\ref{cor:eliminate} with $U = V(G)$ turns it into a strong trace (resp.\ a $d$-stable trace, as $\delta(G) > d$). Hence assume that $E \neq \emptyset$ and $S \neq \emptyset$. By Theorem~\ref{thm:max_flow}, condition~(i) gives an Euler tour of each connected component $B$ of $G_S$, traversing every arc of $B$ in its prescribed direction; traversing it twice gives a closed walk $W_B$ which visits every vertex of $B$ and traverses every edge and arc of $B$ twice in the same direction ($W_B$ need not be a strong trace of $B$). By (ii) and Corollary~\ref{cor:multigraph_tree_trace} (resp.\ Corollary~\ref{cor:multigraph_tree_trace_d}), $G/S$ admits an antiparallel double trace $W^\ast$ whose nontrivial repetitions (resp.\ nontrivial repetitions of order at most $d$) occur only at $S$-vertices. Lemma~\ref{lem:assembly}, applied with $S$ in place of $E'$ and the walks $W_B$, yields a double trace $W$ of $G$ in which every edge of $E'$ is parallel, every edge of $E$ is antiparallel, every arc is traversed twice in its prescribed direction, and whose local transition multigraph at every vertex which is not an $S$-vertex is that of $W^\ast$. Every $S$-vertex is incident with an edge of $E'$ or an arc, hence with an element traversed twice in the same direction, so Corollary~\ref{cor:eliminate} with $U$ the set of $S$-vertices yields a double trace $\widetilde{W}$ with the same traversal directions, without nontrivial repetitions at $S$-vertices, and with the local transition multigraphs of $W^\ast$ elsewhere. Thus $\widetilde{W}$ is an $E$-restricted strong trace (resp.\ every nontrivial repetition of $\widetilde{W}$ has order greater than $d$, and $\delta(G) > d$ excludes trivial repetitions $N(v)$ with $d(v) \leq d$, so $\widetilde{W}$ is an $E$-restricted $d$-stable trace) in which every arc is traversed twice in its prescribed direction.
\end{proof}

\section{Algorithmic aspects and open problems}

The two operations underlying our constructions---the merging of Lemma~\ref{lemma:joining_cycles} and the repetition elimination of Lemma~\ref{lemma:vertex_figures_cycles}---generalize, respectively, the {\em crossing} and {\em $3$-interchange} operations that Benevant~L\'opez and Soler~Fern\'andez~\cite{ben-1998}, building on Thomassen's work~\cite{th-1990}, used to turn a suitable spanning tree into an antiparallel $1$-stable trace in polynomial time. Consequently, the trace-construction half of our characterizations is efficient once a witnessing spanning tree is given.

\begin{proposition}
\label{prp:construction}
Let $G$ be a connected graph with $m = |E(G)|$ edges and $E \subseteq E(G)$ such that condition~{\rm (i)} of Theorem~\ref{thm:BR_strong} holds, and let $T$ be a spanning tree of $G/E'$ witnessing condition~{\rm (ii)} of Theorem~\ref{thm:BR_strong} {\rm (}respectively, of Theorem~\ref{thm:BR_d-stable}, where in addition $\delta(G) > d${\rm )}. Then an $E$-restricted strong trace {\rm (}respectively, an $E$-restricted $d$-stable trace{\rm )} of $G$ can be constructed from $T$ in $O(m^2)$ time, in addition to the time needed for one application of Theorem~\ref{thm:anti_strong} to the graph $G_\xi$ arising in the proof of Lemma~\ref{lemma:tree_with_v_vertices} {\rm (}respectively, of Lemma~\ref{lemma:tree_with_v_vertices_d}{\rm )} and, in the $d$-stable case, for the balanced splittings of Lemma~\ref{lem:balanced_split}. Both of these steps are constructive: the proof of Theorem~\ref{thm:anti_strong} in~\cite{fi-2014} rests on Xuong's theory~\cite{xu-1979}, for which the algorithm of Furst, Gross and McGeoch~\cite{fur-1988} constructs a single-face embedding from a suitable spanning tree, and Lemma~\ref{lem:balanced_split} is obtained in~\cite{rus-2017} from an explicit splitting construction.
\end{proposition}

\begin{proof}
Checking that $G[E']$ is even, forming $G/E'$ and the simple graph $G'$ associated with it in Lemma~\ref{lem:multigraph_strong}, performing the bridge splittings in the proof of Lemma~\ref{lemma:tree_with_v_vertices} (each odd co-tree component is processed once, in time linear in its size), and lifting walks between $G_\xi$, $G'$, $G/E'$ and $G$ all take $O(m)$ time. Transferring $T$ to $G'$ by the proof of Lemma~\ref{lem:multigraph_join_multiple_vertex_trees} processes at most $m$ gadgets, each requiring a connectivity computation in the current co-tree, hence $O(m^2)$ time in total. We represent a closed walk by a circular doubly linked list of edge traversals and store, for every vertex $v$, its visits (pairs of consecutive traversals). Euler tours of the components of $G[E']$ are found in $O(m)$ time. Merging two closed walks at a common vertex (Lemma~\ref{lemma:joining_cycles}) relinks a constant number of pointers; the assembly of Lemma~\ref{lem:assembly} performs $O(m)$ merges after $O(m)$ preprocessing (pairing the segments at boundary vertices, constructing and traversing a spanning tree of the graph $\mathcal{I}$). An application of Lemma~\ref{lemma:vertex_figures_cycles} at a vertex $v$ also relinks a constant number of pointers once the cyclic order along the walk of the three visits involved (the two traversals of $e$ and the $R_2$-visit) is known; determining this order takes $O(m)$ time by traversing the walk. To eliminate the repetitions at a vertex $v \in U$ (Corollary~\ref{cor:eliminate}), we compute the components of $H_v$ once, in $O(d(v))$ time, and then repeatedly apply Lemma~\ref{lemma:vertex_figures_cycles} with $R_1$ the component containing a neighbor $u$ joined to $v$ by a parallel edge $e$, $R_2 = N(v) \setminus R_1$, and the $R_2$-visit chosen in one further component $K$ of $H_v$; each application merges $R_1$ with $K$. Since each application decreases the number of components of $H_v$ by one, at most $d(v) - 1$ applications are needed at $v$, and summing over all vertices gives at most $2m$ applications, hence $O(m^2)$ time. The same accounting applies to the $d$-stable case.
\end{proof}

\medskip
Whether the spanning-tree conditions of Theorems~\ref{thm:BR_strong} and~\ref{thm:BR_d-stable} can themselves be {\em decided} efficiently is a separate question. In the classical case $E' = \emptyset$ and $d = 1$ the condition asks for a spanning tree whose odd co-tree components all contain a vertex of degree at least $4$; without the degree exemption this is the question whether the deficiency $\xi(G)$ vanishes, which is decidable in polynomial time through Xuong's theory~\cite{xu-1979, xu2-1979} and the reduction of Furst, Gross and McGeoch~\cite{fur-1988} to the matroid parity problem solved by Gabow and Stallmann~\cite{gab-1986}. Benevant~L\'opez and Soler~Fern\'andez~\cite{ben-1998} carried out exactly this reduction for the degree-$\geq 4$ exemption. Our conditions differ from theirs only in the {\em exemption rule}: an odd co-tree component is permitted when it contains an $E'$-vertex or a vertex of degree at least $2d + 2$, in place of a vertex of degree at least $4$. This suggests the following.

\begin{conjecture}
\label{conj:tractable}
For every fixed integer $d \geq 1$, deciding whether a connected graph $G$ and a set $E \subseteq E(G)$ satisfy the conditions of Theorem~\ref{thm:BR_strong} {\rm (}respectively, Theorem~\ref{thm:BR_d-stable}{\rm )} can be done in polynomial time.
\end{conjecture}

A natural approach is a reduction to a (co-)graphic matroid parity instance in which parity pairs are attached only to those co-tree components that contain no exempt vertex: the components that already contain an exempt vertex are left unconstrained, so that only the remaining components must be forced to be even. Verifying that this construction is correct---and that it is compatible with the contraction $G/E'$, whose odd components may merge or split exempt vertices---is the one missing ingredient, which we leave to future work.

By way of quantification, the augmenting path algorithm of Gabow and Stallmann~\cite{gab-1986} solves linear matroid parity in $O(mn^3)$ time and graphic matroid parity in $O(mn^2)$ time, where $m$ is the number of elements and $n$ the rank, and for graphic matroids the same authors give an $O(mn \log^6 n)$ algorithm in~\cite{gab-1985}. Combined with the construction of Proposition~\ref{prp:construction}, a positive resolution of Conjecture~\ref{conj:tractable} would yield a polynomial-time procedure that both decides the existence of, and constructs, an $E$-restricted strong trace (or $d$-stable trace)---which is of practical interest for researchers implementing nanostructure design tools.

\section*{Acknowledgments}

Dan Archdeacon passed away in February 2015. While the core of this research was developed during his lifetime, the final manuscript was prepared for publication only recently. We dedicate this work to his memory and to his foundational contributions to the field.



\begin{thebibliography}{99}

\bibitem{ba-1979}
  V.~Batagelj, T.~Pisanski,
  \emph{On partially directed eulerian multigraphs},
  Publ. Inst. Math. (Beograd) (N.S.) 25(39) (1979), 16--24.

\bibitem{ben-1998}
  E.~Benevant~L\' opez, D.~Soler~Fern\' andez,
  \emph{Searching for a strong double tracing in a graph},
  TOP 6 (1998), 123--138.

\bibitem{eg-1984}
  R.~B.~Eggleton, D.~K.~Skilton,
  \emph{Double tracings of graphs},
  Ars Combin. 17A (1984), 307--323.

\bibitem{fi-2014}  
  G.~Fijav\v z, T.~Pisanski, J.~Rus, 
  \emph{Strong traces model of self-assembly polypeptide structures},
  MATCH Commun. Math. Comput. Chem. 71 (2014), 199--212.

\bibitem{fl-1983}  
  H.~Fleischner,  
  \emph{Eulerian graphs}, 
  Selected Topics in Graph Theory 2 (1983),  Academic Press, London-New York, 17--53.

\bibitem{fl-1990}  
  H.~Fleischner, 
  \emph{Eulerian Graphs and Related Topics. Part 1. Vol. 1.}, 
  North-Holland, Amsterdam, 1990.
  
\bibitem{fl-1991}  
  H.~Fleischner, 
  \emph{Eulerian Graphs and Related Topics. Part 1. Vol. 2.}, 
  North-Holland, Amsterdam, 1991. 
  
\bibitem{fo-1962}  
  L.~R.~Ford~Jr., D.~R.~Fulkerson, 
  \emph{Flows in Networks}, 
  Princeton Univ. Press, Princeton, New Jersey, 1962.

\bibitem{fur-1988}
  M.~L.~Furst, J.~L.~Gross, L.~A.~McGeoch,
  \emph{Finding a maximum-genus graph imbedding},
  J. Assoc. Comput. Mach. 35 (1988), 523--534.

\bibitem{gab-1985}
  H.~N.~Gabow, M.~Stallmann,
  \emph{Efficient algorithms for graphic matroid intersection and parity},
  in: Automata, Languages and Programming (ICALP 1985), Lecture Notes in Comput. Sci. 194, Springer, Berlin, 1985, 210--220.

\bibitem{gab-1986}
  H.~N.~Gabow, M.~Stallmann,
  \emph{An augmenting path algorithm for linear matroid parity},
  Combinatorica 6 (1986), 123--150.

\bibitem{gr-2013}
  H.~Gradi\v sar, S.~Bo\v zi\v c, T.~Doles, D.~Vengust, I.~Hafner~Bratkovi\v c, A.~Mertelj, B.~Webb, A.~\v Sali, S.~Klav\v zar, R.~Jerala,
  \emph{Design of a single-chain polypeptide tetrahedron assembled from coiled-coil segments},
  Nature Chemical Biology 9 (2013), 362--366.
  
\bibitem{kl-2013}
  S.~Klav\v zar, J.~Rus,
  \emph{Stable traces as a model for self-assembly of polypeptide nanoscale polyhedrons}, 
  MATCH Commun. Math. Comput. Chem. 70 (2013), 317--330.

\bibitem{ko-2015}
  V.~Ko\v car, S.~Bo\v zi\v c~Abram, T.~Doles, N.~Ba\v si\' c, H.~Gradi\v sar, T.~Pisanski, R.~Jerala,
  \emph{TOPOFOLD, the designed modular biomolecular folds: polypeptide-based molecular origami nanostructures following the footsteps of DNA},
  WIREs Nanomed. Nanobiotechnol. 7 (2015), 218--237.

\bibitem{lj-2017}
  A.~Ljube\v ti\v c, F.~Lapenta, H.~Gradi\v sar, I.~Drobnak, J.~Aupi\v c, \v Z.~Strm\v sek, D.~Lain\v s\v cek, I.~Hafner-Bratkovi\v c, A.~Majerle, N.~Krivec, M.~Ben\v cina, T.~Pisanski, T.~\'Cirkovi\'c Veli\v ckovi\'c, A.~Round, J.~M.~Carazo, R.~Melero, R.~Jerala,
  \emph{Design of coiled-coil protein-origami cages that self-assemble in vitro and in vivo},
  Nature Biotechnology 35 (2017), 1094--1101.

\bibitem{moh-2001}
  B.~Mohar, C.~Thomassen,
  \emph{Graphs on Surfaces},
  The Johns Hopkins University Press, 2001. 
  
\bibitem{ore-1951}
  O.~Ore,
  \emph{A problem regarding the tracing of graphs},
  Elemente der Math. 6 (1951), 49--53.
 
\bibitem{rus-2017} 
  J.~Rus,
  \emph{Antiparallel $d$-stable traces and a stronger version of Ore problem},
  J. Math. Biol. 75 (2017), 109--127.

\bibitem{sa-1977}
  G.~Sabidussi,
  \emph{Tracing graphs without backtracking},
  Operations Research Verfahren XXV, Symp. Heidelberg, Teil 1 (1977), 314--332.

\bibitem{th-1990}
  C.~Thomassen,
  \emph{Bidirectional retracting-free double tracings and upper embeddability of graphs},
  J. Combin. Theory Ser. B 50 (1990), 198--207. 

\bibitem{tro-1966}
  D.~J.~Troy,
  \emph{On traversing graphs},
  Amer. Math. Monthly 73 (1966), 497--499.

\bibitem{ve-1975}
  P.~D.~Vestergaard, 
  \emph{Doubly traversed euler circuits},
  Arch. Math. 26 (1975), 222--224. 

\bibitem{wa-1970}
  K.~Wagner, 
  \emph{Graphentheorie},
  BI-Hochschultaschenb\"ucher, Bd. 248, Bibliograph. Inst. AG, Mannheim, 1970.

\bibitem{wan-2017}
  J.~Wang, G.~Hu, M.~Ji,
  \emph{Almost parallel strong trace model of self-assembly polypeptide nanostructure},
  MATCH Commun. Math. Comput. Chem. 77 (2017), 783--798.
  
\bibitem{wan2-2017}
  J.~Wang, X.~Jin, M.~Ji,
  \emph{On the existence of $F$-strong trace of a graph when $F$ induces a forest},
  MATCH Commun. Math. Comput. Chem. 77 (2017), 799--812.

\bibitem{we-1996}
  D.~B.~West, 
  \emph{Introduction to Graph Theory},
  Prentice Hall, Upper Saddle River, 1996. 
  
\bibitem{xu-1979}
  N.~H.~Xuong,
  \emph{How to determine the maximum genus of a graph}, 
  J. Combin. Theory Ser. B 26 (1979), 217--225. 

\bibitem{xu2-1979}
  N.~H.~Xuong,
  \emph{Upper-embeddable graphs and related topics}, 
  J. Combin. Theory Ser. B 26 (1979), 226--232. 
   
\end{thebibliography}
\end{document}